\documentclass{article}
\usepackage{setspace}
\usepackage[english]{babel}
\usepackage{mathrsfs}
\usepackage{tikz}
\usepackage[english]{babel}
\usepackage[left=1 in, right=1 in, top=1.5 in, bottom=1.5 in]{geometry}
\usepackage{amsmath}
\usepackage{amsthm}
\usepackage{amssymb}
\usepackage{mathtools}
\usepackage[utf8]{inputenc}
\usepackage{color}
\usepackage{amsthm}
\usepackage{tikz}
\usetikzlibrary{graphs,arrows.meta}
\usepackage[shortlabels]{enumitem}
\usepackage{graphicx}
\usepackage{caption}
\usepackage{subcaption}
\usepackage{hyperref}
\usepackage{tikz}
\usepackage{ytableau}
\newtheorem{theorem}{Theorem}[section]
\newtheorem{lemma}{Lemma}[section]
\newtheorem{proposition}{Proposition}[section]
\newtheorem{corollary}{Corollary}[theorem]  
\newtheorem{definition}{Definition}[section]
\newtheorem{example}{Example}[section]
\newtheorem{remark}{Remark}[section]
\usepackage{forest}
\usepackage{array}
\usepackage{geometry}
\newtheorem{notation}[theorem]{Notation}
\usepackage{forest}
\newcommand{\dsc}{\mathrel{|}}   
\newcommand{\asc}{,\,}           

\newcommand{\N}{\mathbb{N}}

\newcommand{\vanish}[1]{}
\newcommand{\wmax}[2][S]{\overline{w}^{\,#1}_{#2}}
\newcommand{\wmin}[2][S]{\underline{w}^{\,#1}_{#2}}
\usepackage{amsmath}
\usepackage{graphicx}
\usepackage{lineno}

\title{Polynomiality of Subdimensions of Diagonal Harmonics and a Sharp Stability Bound}
\author{
  Xinxuan Wang \\
  Department of Mathematics \\
  University of Pennsylvania \\
}

\begin{document}
\maketitle
\renewcommand{\thefootnote}{}%

\begin{abstract}
A sequence of $S_n$-representations $\{V_n\}_{n \ge 1}$ is
\emph{representation stable} if, writing
$V(\lambda) = (n-|\lambda|, \lambda_1, \lambda_2, \dots)$ for each
partition $\lambda$, the multiplicity of the irreducible indexed by
$V(\lambda)$ in $V_n$ is eventually independent of $n$. In particular, Church, Ellenberg and Farb \cite{Church_2015} found that if we fix $a$ and $b$, then the space of diagonal harmonics $DH_n^{a,b}$ exhibits this behavior, and its dimension stabilizes to a polynomial in $n$ eventually. Building on this result, we use the Schedules Formula by Haglund and Loehr \cite{HAGLUND2005189} to get an explicit combinatorial polynomial for the dimension of the bigraded spaces $DH_n^{a,b}$. This derivation not only yields the dimension formula but also produces a new stability bound of \( a + b \) which is sharp, and determines the exact degree of the dimension polynomial, which is also \( a + b \).

\end{abstract}
\section{Introduction}
\subsection{Overview}
In this paper, we begin by introducing the Diagonal Harmonics, the central object of our study, in Section~\ref{sec:diagharmonics}. We then introduce Schedules Formula, the combinatorial formula of Haglund and Loehr \cite{HAGLUND2005189} that we use to derive the formula in Section~\ref{sec:schedules}. Next, we give a brief overview of representation stability, a phenomenon discovered by Church and Farb \cite{Church_2013}, which serves as a key motivation for this work in Section~\ref{sec:respstab}. In particular, it suggests that the dimension of the bigraded subspace of $DH_n$ is given by a polynomial in $n$. \\
Section~\ref{sec:formula} contains our main result: an explicit polynomial expression for the dimensions of the bigraded components. Finally, in Section~\ref{sec:sharp} we prove that $a+b$ is a sharp bound for stabilization, and that the degree of the dimension polynomial of $DH_n^{a,b}$ is $a + b$.
\subsection{Diagonal Harmonics}\label{sec:diagharmonics}
\begin{definition}
    Given a $f \in \mathbb{C}[x_1,\dots,x_n]$, $S_n$ acts on $f$ by permuting the indices of the variables, i.e.
\begin{equation*}
    \sigma \cdot f(x_1,...,x_n) = f(x_{\sigma_1},...,x_{\sigma_n})
\end{equation*}
    The \textbf{Coinvariant Ring of $S_n$} is 
            \begin{equation}
                R_n = \frac{\mathbb{C} [x_1,\dots,x_n] }{I(x_1,\dots,x_n)^+}
            \end{equation}
            where $x_1,\dots,x_n = \{x_1,x_2,...,x_n\}$, and 
$I(x_1,\dots,x_n)^+ = \{f \in \mathbb{C}[x_1,\dots,x_n] | \sigma \cdot f = f \text{ and } f \text{ is not a constant}\}$.
\end{definition}
In addition to its rich ring structure, the coinvariant ring $R_n$ is
also isomorphic to the singular cohomology ring of the complete flag
variety $\mathrm{Fl}_n$. In other words,
\[
R_n \cong H^*(\mathrm{Fl}_n;\mathbb{C}).
\]

It is also known to be isomorphic to the space of harmonics $H_n$ as $S_n$-module \cite{haiman1994conjectures}, which we define below.
 \begin{definition}
     \textbf{The Harmonics} are
            \begin{equation}
                H_n = \{ f(x_1,\dots,x_n) \in \mathbb{C} [x_1,\dots,x_n] :\sum\limits_{i=1}^n \partial_{x_i}^k f(x_1,\dots,x_n) = 0 \text{ for all k} > 0\}
            \end{equation} 
 \end{definition}
Let $R_n^d$ denote the degree-$d$ graded component of
$R_n$, consisting of the elements of $R_n$ represented by homogeneous
polynomials of degree $d$ in $R_n$.
Artin~\cite{artin1944galois} found a basis for $R_n$, which makes it
easy to determine the dimension of $R_n$ and its graded components $R_n^d$.
\begin{theorem}
    $\{x_1^{\alpha_1} \cdots x_n^{\alpha_n} : 0 \le \alpha_i < i \}$ form a basis for $R_n$
\end{theorem}
\begin{corollary}
    $dim(R_n)$ = n!
\end{corollary}
The Corollary follows from directly computing the number of basis elements. 
\begin{definition}
    \textbf{The Hilbert Series of $R_n$} is defined to be 
    \begin{align}
        Hilb(R_n) & = \sum\limits_{d \ge 0} q^d \cdot dim(R_n^d)
    \end{align}
\end{definition}
\begin{corollary} Recall that a composition of $n$ is a sequence $\alpha=(\alpha_1,\ldots,\alpha_k)$ of positive integers such that $\alpha_1+\cdots+\alpha_k=n$. We write $\alpha\models n$ to indicate that $\alpha$ is a composition of $n$. Then by Artin's basis, the Hilbert Series of $R_n$ is 
    \begin{align}
        Hilb(R_n) &= \sum\limits_{\substack{\alpha \models n \\ 0 \le \alpha_i < i }} q^{|\alpha|} \\
        &= 1 \cdot (1+q) \cdot (1+q+q^2) ... (1+q+...q^{n-1}) \label{eq5}
    \end{align}
\end{corollary}
If we introduce a new notation called the $q$-integers where $[k]_q = 1 + q + q^2 + ... + q^{k-1}$, then we can write (\ref{eq5}) as $[n]_q! = [n]_q[n-1]_q...[1]_q$, and the dimensions of the bigraded component $R_n^d$ would be the coefficient of $q^d$ in $[n]_q!$.\\
Beyond the classical coinvariant ring $R_n$, there are several families
of generalized coinvariant rings that have been studied in algebraic and
enumerative combinatorics; see, for example,
\cite{HaglundRhoadesShimozono2018,RhoadesWilson2017,RhoadesYuZhao2020}.
The diagonal coinvariant ring $DR_n$ is one such example. \\
Given a polynomial $f(x_1,\dots,x_n,y_1,\dots,y_n) \in \mathbb{C}[x_1,\dots,x_n,y_1,\dots,y_n]$, $S_n$ acts diagonally by permuting the $x$ and $y$ variables by
            \begin{equation}
                \sigma\cdot f(x_1,...,x_n,y_1,...,y_n) = f(x_{\sigma_1},...,x_{\sigma_n},y_{\sigma_1},...,y_{\sigma_n})
            \end{equation}
\begin{definition}
\textbf{The Diagonal Coinvariant Ring of $S_n$} is 
            \begin{equation}
                DR_n = \frac{\mathbb{C}[x_1,\dots,x_n,y_1,\dots,y_n]}{I(x_1,\dots,x_n,y_1,\dots,y_n)^+}
            \end{equation}
            where $I(x_1,\dots,x_n,y_1,\dots,y_n)^+$ is the ideal generated by all polynomial\\ $f(x_1,\dots,x_n,y_1,\dots,y_n) \in \mathbb{C}[x_1,\dots,x_n,y_1,\dots,y_n]$ that are invariant under the diagonal action of $S_n$ without constant term.    
\end{definition}
$DR_n$ also has interesting geometric connections. Carlsson and
Oblomkov~\cite{CarlssonOblomkov2018} related $DR_n$ to type $A$ affine
Springer fibers using the Lusztig--Smelt paving of these varieties,
obtaining an independent proof of the dimension formula
$\dim(DR_n)=(n+1)^{n-1}$ by Haiman \cite{haiman1994conjectures} as well as an explicit monomial basis of
$DR_n$. The number $(n+1)^{n-1}$ is also the classical enumeration of
parking functions of length $n$~\cite{Stanley1996}. By Cayley's formula,
the same number counts labeled trees on $n+1$ vertices, and parking
functions are known to be in bijection with such labeled trees
\cite{IrvingRattan2021}. In our setting, the parking-function
interpretation is the one that will be relevant below. \\
Similar to $R_n$ being isomorphic to $H_n$ as an $S_n$-module, $DR_n$ is also isomorphic
to $DH_n$ as an $S_n$-module~\cite{haiman1994conjectures}, which we
define below.
\begin{definition}
\textbf{The Space of Diagonal Harmonics} is
\begin{align*}
DH_n
= \biggl\{ f(x_1,\dots,x_n,y_1,\dots,y_n)
&\in \mathbb{C}[x_1,\dots,x_n,y_1,\dots,y_n] :\\
&\sum_{i=1}^n
\partial_{x_i}^r\partial_{y_i}^s
f(x_1,\dots,x_n,y_1,\dots,y_n)=0
\text{ for all } r,s\geq 0 \text{ with } r+s>0
\biggr\}.
\end{align*}
\end{definition}
The dimension of $DR_n$ becomes much harder to track when we add another set of variables, and the original and first proof of its dimension was by Haiman using geometry of Hilbert Schemes. 
\begin{theorem}
    \cite{Haiman_2002}
    \begin{align}
    dim(DR_n) = (n+1)^{n-1}
\end{align}
\end{theorem}
The derivation of the bigraded Hilbert series is also nontrivial, as it
remained a conjecture for quite some time until Carlsson and
Mellit~\cite{carlsson2018proofshuffleconjecture} proved the Shuffle
Theorem. We first describe the bigrading of $DR_n$. For $a,b\geq 0$,
let $DR_n^{a,b}$ denote the bihomogeneous component of $DR_n$ consisting
of polynomials that are homogeneous of degree $a$ in the variables
$x_1,\dots,x_n$ and homogeneous of degree $b$ in the variables
$y_1,\dots,y_n$. We then define the Hilbert series of $DR_n$ as follows.
\begin{definition}
    The \textbf{Hilbert Series} of $DR_n$ is
\begin{align}
     Hilb(DR_n) &= \sum\limits_{a,b \ge 0} q^at^b \cdot dim\big(DR_n^{a,b}\big)
\end{align}
\end{definition}
The Shuffle Theorem gives a combinatorial formula for the bigraded
Hilbert series of $DR_n$ in terms of parking functions and their
statistics $\operatorname{area}$ and $\operatorname{dinv}$
\cite{carlsson2018proofshuffleconjecture}. This formula was subsequently
reformulated by Haglund and Loehr \cite{HAGLUND2005189} as the Schedules Formula, which is the
formulation we use in this paper.
\section{The Schedules Formula} \label{sec:schedules}
Before we present the Schedules Formula, we need to introduce two permutation statistics, $maj(\sigma)$ and $w_i(\sigma).$
\begin{definition}
    The \textbf{major index} statistic $maj(\sigma)$ of a permutation $\sigma \in S_n$ is 
    \begin{equation}
        maj(\sigma) = \sum_{i \text{ s.t. }\sigma_i > \sigma_{i+1}} i
    \end{equation}
    Furthermore, $i$ is called a descent if $\sigma_i > \sigma_{i+1}$.
\end{definition}
\begin{definition}
For $\sigma \in S_n$, let the 1-st run be defined to be the first increasing subsequence, i-th run be defined to be the i-th increasing subsequence.
   \begin{itemize}
    \item $w_i(\sigma) = $\#(entries $\sigma_j$ in the same run as $\sigma_i$ and $\sigma_j > \sigma_i \bigr) + \#\bigl($entries $\sigma_k$ in the next run of $\sigma_i$ and $\sigma_k < \sigma_i \bigr)$ 
    \end{itemize}
    While calculating $w_i$, we adjoint 0 at the end of $\sigma$. Denote the sequence $\big(w_1(\sigma), \dots, w_n(\sigma)\big)$ as $wseq(\sigma)$ and call it the \textbf{w-sequence of $\sigma$}.
\end{definition}
\begin{theorem} (The Schedules Formula) \cite{HAGLUND2005189}
    \begin{align} \label{eq:schedules}
    Hilb(DR_n) =  \sum\limits_{\sigma \in S_n} t^{\text{maj}(\sigma)} \prod\limits_{i=1}^n [w_i(\sigma)]_q 
\end{align}
\end{theorem}
\begin{example}
   \begin{align}
    \sigma &= (4,2,5,1,3,8,6,7,9) \in S_9\\ 
    S &= \{1,3,6\}\\
    maj(\sigma) &= 1+3+6 = 10\\
    wseq(\sigma) &= (1, 2,2,2,1,2,3,2,1)
\end{align}
Putting everything together, $\sigma$ gives $t^{10} [2]^4 \cdot [3] \cdot [2] \cdot [1]$.
\end{example}
Introduce the notation $[q^k](f)$ for the coefficient of $q^k$ in $f$,
then
\(
\dim(DR_n^{a,b}) = [q^a t^b]\bigl(\mathrm{Hilb}(DR_n)\bigr).
\)
The Schedules Formula provides the combinatorial framework that we use
to derive an explicit formula for $\dim(DR_n^{a,b})$. The motivation for
seeking such a formula comes from representation stability: as we will
see in the next section, representation stability implies that
$\dim(DR_n^{a,b})$ is eventually polynomial in $n$. Thus, while
representation stability tells us that such a polynomial must exist, the
Schedules Formula allows us to determine it explicitly.
\section{Representation Stability} \label{sec:respstab}
The irreducible representations of $S_n$ are indexed by partitions
$\lambda = (\lambda_1, \lambda_2, \dots, \lambda_r) \vdash n$, where
$\lambda_1 \geq \lambda_2 \geq \dots \geq \lambda_r > 0$ and
$\sum_i \lambda_i = n$~\cite{sagan2001symmetric}. Following notations of Church and Farb, let $V(\lambda) = V(\lambda_1,...,\lambda_r)$ be the irreducible $S_n$-representation corresponding to the partition $\bigl((n-\sum_{i=1}^r\lambda_i), \lambda_1,...,\lambda_r\bigr),$ where $\lambda_1 \ge \lambda_2 \ge ... \ge \lambda_r$. Note that this is a partition only when $n - (\sum_{i=1}^r\lambda_i) \ge \lambda_1 $. Church and Farb \cite{Church_2013} studied a phenomenon called representation stability, which is the foundation of our result. Independently, Bergeron \cite{BERGERON201397} developed a notion of representation stability using symmetric function theories, and more recently Borisov \cite{borisov2025monomialstabilityfrobeniusimages} obtained new results in representation stability also through the lens of symmetric function theories. We introduce the basic definitions and results from representation stability that will be used in this paper below.
\begin{definition} \cite{Church_2015}
    A sequence of representations $V_n$ are \textbf{representation multiplicity stable} if the decomposition of $V_n$ into $S_n$-representations as 
    \begin{align}
        V_n = \bigoplus_\lambda c_{\lambda, n} V(\lambda)
    \end{align}
    stabilizes, i.e. the coefficient $c_{\lambda ,n}$ are eventually independent of $n$.
\end{definition}
A more refined notion of representation stability, developed through the theory of FI-modules, can be found in \cite{Church_2015}. Below is an example illustrating representation stability.
\begin{example}
    The $n$-th Configuration Space of $\mathbb{C}$ is      \begin{equation}
                    Conf_n(\mathbb{C}) = \{(z_1,...,z_n) \in \mathbb{C}^n | \forall i < j, z_i \ne z_j \}
                \end{equation}
    There is an $S_n$ action on $Conf_n(\mathbb{C})$, which induces an
$S_n$ action on its cohomology, making $H^i(Conf_n(\mathbb{C}))$ an
$S_n$-representation. Moreover, Church and Farb
shows that the sequence $\{H^i(Conf_n(\mathbb{C}))\}_{n\geq 1}$ is
representation stable. Below is the calculation of
$H^2(Conf_n(\mathbb{C}))$ due to Church and Farb.
    \begin{align*}
                        H^2(Conf_5(\mathbb{C})) &= V(1)^{\oplus 2} \oplus V(1,1)^{\oplus 2} \oplus V(2)^{\oplus 2} \oplus V(2,1)\\
                        H^2(Conf_6(\mathbb{C})) &= V(1)^{\oplus 2} \oplus V(1,1)^{\oplus 2} \oplus V(2)^{\oplus 2} \oplus V(2,1)^{\oplus 2} \oplus V(3)\\
                        H^2(Conf_n(\mathbb{C})) &= V(1)^{\oplus 2} \oplus V(1,1)^{\oplus 2} \oplus V(2)^{\oplus 2} \oplus V(2,1)^{\oplus 2} \oplus V(3) \oplus V(3,1)
                    \end{align*} for $n \ge 7$.
    \end{example}
\begin{theorem} \cite{Church_2015}
    If $V_n$ is representation stable, then $dim(V_n)$ is eventually a polynomial in $n$.
\end{theorem}
\begin{proof}
    This result can be derived using two perspectives; one using the notion of character polynomials, the other using the hook length formula. For the proof using character polynomials, see \cite{Church_2015}. \\
    For $\lambda = (\lambda_1,\lambda_2, \dots, \lambda_r)$, let $\lambda(n) = (n - |\lambda|, \lambda_1, \lambda_2, \dots, \lambda_r)$. Recall that by hook formula, dimension of $V(\lambda_1,...,\lambda_r)$ is given by 
                
                    \begin{align}
                    f^{\lambda(n)} &= \frac{n!}{\prod\limits_{(i,j) \in \lambda(n)} h_{\lambda(n)}(i,j)}\\
                    &= \frac{n!}{\prod\limits_{(1,j) \in \lambda} h_\lambda(i,j)\cdot \prod\limits_{(i, j) \in \lambda, i > 1} h_\lambda(i,j)} \label{eq15} \\
                    &= \frac {n!}{(n-|\lambda| -\lambda _1)! \prod_{j=1}^{\lambda _1} (n-|\lambda|  -j+1 + \lambda  ^{\prime}_j )}
\times \frac {f^{\lambda}}{|\lambda|!}
                \end{align}
    There are $n - |\lambda|$ many terms containing $n$ in the first product $\prod\limits_{(1,j) \in \lambda} h_\lambda(i,j)$. The factors are distinct and cancel distinct terms in $n!$, thus leaving behind a product of $\prod_{k \in S} (n-k)$ for some $S \subseteq \{1,2,...,n\}$, so (\ref{eq15}) is always a polynomial in $n$ of degree $|\lambda|$, given that $n$ is big enough for $\lambda(n)$ to be a partition, i.e. when $n > |\lambda| + \lambda_1$.
\end{proof}
\begin{theorem} \cite{Church_2015}
    If $a$ and $b$ are fixed, then $DR_n^{a,b}$ are representation stable, and $dim(DR_n^{a,b})$ is eventually a polynomial in $n$.
\end{theorem}
While the theory of representation stability ensures that $\dim(DR_n^{a,b})$ is eventually polynomial in $n$, it does not specify the polynomial explicitly. Section~\ref{sec:formula} addresses this gap by providing a concrete formula using the Schedules Formula. From the general framework of representation stability, it follows that the degree of this polynomial is at most $a + b$. 
Bahran \cite{BAHRAN2024429} provided a stability range using FI-module theory to be $n \ge 2(a+b)-1$.
However the theory does not yield a sharp stabilization bound, and no general method is known for proving sharpness using representation stability theories. In Section~\ref{sec:sharp}, we establish that $a + b$ is a sharp bound, and that the degree of the polynomial is precisely $a+b$ — both immediate consequences of the structure of our explicit polynomial.

\section{Explicit Formula for $\dim(DR_n^{a,b})$} \label{sec:formula}
\subsection{Formula and Preliminary Setup}
\begin{definition} 
Let $n,b\in\mathbb{N}$, $\sigma\in S_n$, $\tau\in\mathbb{N}^n$, and let
$S=\{s_1,\dots,s_d\}\subseteq\{1,2,\dots,n-1\}$ with
$s_1<s_2<\dots<s_d$ and $s_1+s_2+\dots+s_d=b$.
Let $U\subseteq\{1,2,\dots,s_d\}$. Define
    \begin{align*}
                    W(\tau) &= \{\sigma \in S_n| \,wseq(\sigma) = \tau\}\\
                    D_S &= \{\sigma \in S_n | \,Desc(\sigma) = S\}\\
                        W(\tau,U) &= \{\sigma \in S_n | \, w_i(\sigma) = \tau_i \text{ for } i\not\in U, w_i(\sigma) \ge \tau_i \text{ for } i \in U\}
                    \end{align*} 
\end{definition}
\begin{theorem}\label{thm:main_theorem}
    (W.) 
For $a\in\mathbb{Z}_{\ge0}$ and $b\in\mathbb{N}$, and for $n$ sufficiently
large relative to $a$ and $b$, the dimension of $DH_n^{a,b}$ is
\begin{small}
                        \begin{align} \label{formula2}
                    \sum\limits_{\substack{S \subseteq \{1,2,...,b\} \\ \text{ }s_1+s_2+...+s_d = b}}
   \sum_{U \subseteq \{1,2,...,s_d\}}
   \sum_{k=0}^a
   \sum\limits_{\substack{
    \tau_l = k+1 \text{ for } l \in U \\
    \tau_l \le k \text{ for } l \notin U \\
    \tau \text{ permissible}
}}
   \Bigg(\bigg|D_S \bigcap W(\tau;U)\bigg|\Bigg) \cdot 
   \Bigg([q^k] \bigg(\prod_{i=1}^{s_d} [\tau_i]_q\bigg)\Bigg)  \cdot 
   \Bigg([q^{a-k}]\bigg([n-s_d]_q!\bigg)\Bigg)
                \end{align}
            \end{small} 
    In particular, $\dim DH_n^{a,b}$ is a polynomial in $n$.
\end{theorem}
Observe that the tail of the $w$-sequence is completely determined by $s_d$. Indeed,
$w_i = n - i + 1$ for every $i > s_d$, so
\begin{equation}\label{eq:tail}
    \prod_{i = s_d+1}^{n} [w_i]_q = [n-s_d]_q! .
\end{equation}
Evaluating the third factor of \eqref{formula2} therefore reduces to
extracting a coefficient of a $q$-factorial, which is the content of the
following classical result.

\begin{theorem}[Knuth's Formula]\label{thm:knuth}
\cite{knuth1997art}
For $k \le n$,
\begin{equation}\label{eq:knuth}
\begin{split}
    [q^k]\left([n]_q!\right) = \binom{n+k-1}{k}
    &+ \sum_{j=1}^{\infty} (-1)^j \binom{n+k-u_j-1}{k-u_j}
    + \sum_{j=1}^{\infty} (-1)^j \binom{n+k-u_j-j-1}{k-u_j-j},
\end{split}
\end{equation}
where $u_j = \frac{j(3j-1)}{2}$ are the pentagonal numbers.
\end{theorem}
Margolius \cite{Margolius} observed that the sums in \eqref{eq:knuth} are
in fact finite, with $\lfloor \sqrt{\tfrac{1}{36}+\tfrac{2k}{3}} \rfloor$
nonzero terms between them. Hence $[q^k]\left([n]_q!\right)$ is a
polynomial in $n$ for $k \le n$.
In the case $b=0$, the coefficient of $q^a t^0$ in the Schedules formula receives contributions only from permutations with $\operatorname{maj}(\sigma)=0$ and thus no descents; the only such permutation is $\sigma = (1,2,\dots,n-1, n)$, which has $ w_i(\sigma)=i$ for all $i$. Hence \( \prod_{i=1}^n w_i(\sigma)=[n]_q!. \) Thus, the case $b=0$ can be treated directly using Knuth's formula.
\begin{theorem}
    For $a \in \mathbb{Z}_{\ge0}$, $n \ge  a$, dimension of $DR_n^{a,0}$ is given by the Knuth Formula above, i.e.
    \begin{align} \label{b0formula}
        dim(DR_n^{a,0}) = [q^a]\bigg([n]_q!\bigg) = \binom{n+a-1}{a} &+ \sum\limits_{j=1}^\infty (-1)^j \binom{n+a-u_j-1}{a-u_j}
                +\sum\limits_{j=1}^\infty (-1)^j \binom{n+a-u_j-j-1}{a-u_j-j}
    \end{align}
    where $u_j = \frac{j(3j-1)}{2}$ are the pentagonal numbers. It is a polynomial in $n$.
\end{theorem}
Now assume $b\geq 1$. Throughout the paper, we refer to $a$ as the
$q$-degree and to $b$ as the $t$-degree. \\
Given $k \le n$ and $\sigma \in S_n$, define a truncation map as follow:
\[
\text{tr}(y) = 
\begin{cases}
    k + 1, & \text{if } y > k + 1, \\
    y, & \text{otherwise}.
\end{cases}
\]
Next observe that 
\begin{equation}
    [q^k]\bigg( \prod\limits_{i=1}^{s_d} [w_i(\sigma)]_q \bigg)
    = [q^k]\bigg( \prod\limits_{i=1}^{s_d} [tr\big(w_i(\sigma)\big)]_q \bigg)
\end{equation}
With this observation in mind, we only need to focus on the counting problem being: Fix a descent set $S$, given a w-sequence $\tau = (\tau_1, \dots,\tau_n)$ with maximum $\tau_i$ being $k+1$, how many $\sigma \in S_n$ has a w-sequence that truncate to $\tau$?\\
This is exactly the first factor of \eqref{formula2}, namely
$\big|D_S \cap W(\tau,U)\big|$. We claim it is a polynomial in $n$;
this is stated as Theorem~\ref{maincount} at the end of this section
and proved in Section~\ref{sec:mainproof}. Granting it, together with
the fact that permissibility depends only on the descent set $S$ and
guarantees 
\begin{align}
\bigcup_{\substack{\sigma\in S_n\\ \operatorname{Des}(\sigma)=S}}
\bigl(w_1(\sigma),\ldots,w_n(\sigma)\bigr)
=
\bigcup_{\tau\text{ permissible for }S}
(\tau_1,\ldots,\tau_n).
    \end{align}
which will be proved in Theorem~\ref{thm:permissareallwsequence}, we can now prove
Theorem~\ref{thm:main_theorem}.
\begin{proof}[Proof of Theorem \ref{thm:main_theorem}]
    First, notice that the first factor  $\big(\big|D_S \bigcap W(\tau,U)\big|\big)$ is a polynomial in $n$ by Theorem \ref{maincount}, which we assume now and prove later in Section~\ref{sec:mainproof}. The second factor $\Big([q^k] \big(\prod_{i=1}^{s_d} [\tau_i]_q\big)\Big)$ is a constant, and third factor $\Bigl([q^{a-k}]\big([n-s_d]_q!\big)\Bigr)$ is a polynomial in $n$ by Theorem~\ref{thm:knuth} and \cite{Margolius}. \\
    Now we analyze the 4 summations. First, given a fixed $b$ we want to sum over all possible descent set that gives $s_1+s_2+...s_d = b$. Given a descent set $S$, we want to sum over all possible
$\tau \in \mathbb{N}^n$ for which there exists some
$\sigma \in S_n$ with $wseq(\sigma)=\tau$. There are two observations
that allow us to organize this sum so that its indexing set is
independent of $n$. First, for indices $l\in U$, we truncate the
condition $\tau_l\geq k+1$ by replacing it with $\tau_l=k+1$. Thus, the
values of $\tau_l$ that contribute to the sum are bounded independently
of $n$, rather than ranging potentially up to $n$. Second, the entries
of the last run are necessarily
\(
(n-s_d,n-s_d-1,\ldots,2,1).
\)
Although the length of this run is $n-s_d$ and therefore depends on $n$,
its entries are completely determined by $s_d$ and introduce no
additional summation over $n$-dependent choices.
The remaining restriction on $\tau$ comes from permissibility, which we
defer until after this proof. Permissibility depends only on $S$ and ensures that $\tau=(w_1(\sigma),\dots,w_n(\sigma))$ for some $\sigma\in S_n$ by Proposition~\ref{prop:permissible}, as well as the validity of Theorem~\ref{maincount}. Thus, it restricts the possible $\tau$ appearing in the fourth summation. Proposition~\ref{prop:wiispermiss} further shows that every $w$-sequence arising from some $\sigma\in S_n$ is represented by a permissible $\tau$. Moreover, requiring $\tau_l>k$ only when $l\in U$ prevents double counting. Consequently, the fourth summation ranges over a finite index set independent of $n$ that contains exactly the $w$-sequences arising from permutations in $S_n$. This is precisely the content of Theorem~\ref{thm:permissareallwsequence}. Finally, observe that  
    \begin{equation}
        [q^a]\big(\prod\limits_{i=1}^n \big [w_i(\sigma)]_q\big)= \sum_{k=0}^a \bigg([q^k]\big(\prod\limits_{i=1}^{s_d} [w_i(\sigma)]_q\big)\bigg) \cdot \bigg([q^{a-k}]\big(\prod\limits_{i=s_d+1}^n [w_i(\sigma)]_q \big)\bigg)
    \end{equation}
    Therefore we have the third summation. We have established that none of the summations depend on $n$, as they only depend on $a$ and $b$, so we have that (\ref{formula2}) is indeed a polynomial in $n$ that calculates the coefficient of $q^at^b$ in the Schedules Formula.
\end{proof}
We now return to the definition of permissibility and introduce some
additional notation and preliminary results needed to state
Theorem~\ref{maincount}, whose proof will be given in the next section.
\begin{definition}\label{def:bounds}
Let $S=\{s_1<\dots<s_d\}$ be a descent set, and set $s_{d+1}:=n$. For $i\in\{1,\dots,n\}$, let $t$ be the unique index such that $s_t \le i < s_{t+1}$.
Define
\[
  \wmax{i} :=
  \begin{cases}
    (s_{t+1} - i) + (s_{t+2} - s_{t+1}), & t \le d-1, \\
    n - i+1, & t = d,
  \end{cases}
  \qquad
  \wmin{i} :=
  \begin{cases}
    s_{t+1} - i, & i \notin S, \\
    1, & i \in S.
  \end{cases}
\]

\end{definition}

\noindent
The superscript records that these quantities are read off from the
descent set $S$. They bound $w_i$ uniformly on $D_S$:
for every $\sigma \in D_S$ we have
$\wmin{i} \le w_i(\sigma) \le \wmax{i}$ (Proposition~\ref{prop:wiispermiss}).
\begin{definition}
Let a descent set $S$ be given and let $\tau \in \mathbb{N}^{n}$. Write
$R_y = \{\sigma_{s_{y-1}},\, \sigma_{s_{y-1}+1},\, \dots,\, \sigma_{s_{y}}\}$ for the $y$-th run, where we set $s_0=1$ and $s_{d+1}=n$. We call
$\tau$ \textbf{permissible} if
\begin{enumerate}
  \item $\wmin{i} \le \tau_i \le \wmax{i}$ for every $i \in \{1,2,\dots,n-1,n\}$, and
  \item $\tau_j \le \tau_{j+1} + 1$ whenever $j$ and $j+1$ lie in the
        same run $R_y$.
\end{enumerate}
\end{definition}

\begin{definition}\label{def:maxspot}
Let $S$, $\tau$ and $U$ be given with $\tau$ permissible. A position $k$
is a \textbf{maximal spot} if either
\begin{enumerate}
  \item $k = n$ and $\tau_{s_d} < \wmax{s_d}$; or
  \item $k = s_m$ for some $m$, with $s_{m-1} \ne s_m - 1$,
        $\tau_{s_{m-1}} < \wmax{s_{m-1}}$, and either
        $\tau_{s_m} = \wmax{s_m}$ or $k \in U$.
\end{enumerate}

\end{definition}

\begin{proposition}\label{propmaxspot}
Let $\sigma \in D_S \cap W(\tau,U)$. Then the entry $n$ occupies a
maximal spot.
\end{proposition}

\begin{proof}
Since $n$ exceeds every other entry, it can lie only at a descent
position or at position $n$. We rule out the remaining possibilities.

\emph{Position $k = n$.} If $\tau_{s_d} = \wmax{s_d}$ then
$\sigma_{s_d}$ exceeds every entry of the final run, hence exceeds $n$,
which is impossible.

\emph{Position $k = s_m$.} If $s_{m-1} = s_m - 1$ then there is a descent before $\sigma_m = n,$
which is impossible. If
$\tau_{s_{m-1}} = \wmax{s_{m-1}}$ then
$\sigma_{s_{m-1}} > \sigma_{s_m} = n$, also impossible.

\emph{The case $k \in U$.} Since $k \in U$ the constraint at $k$ reads
$w_k(\sigma) \ge \tau_k$ rather than $w_k(\sigma) = \tau_k$, so
$\tau_{s_m}$ does not determine $w_{s_m}(\sigma)$ and the test
$\tau_{s_m} = \wmax{s_m}$ cannot be applied; clause~(2) therefore admits
such a $k$ outright.
\end{proof}
\begin{proposition}
    
\label{prop:maxspotexists}
Let $S$ be given and $\tau$ be permissible. Then at least one
position is a maximal spot.
\end{proposition}

\begin{proof}
If $\tau_{s_d} < \wmax{s_d}$ then $k = n$ is a maximal spot. Otherwise
$\tau_{s_d} = \wmax{s_d}$. Let $v$ be the largest index with
$\tau_{s_{v-1}} < \wmax{s_{v-1}}$ and that $s_{v-1} \ne s_v - 1$, then $s_v$ is a maximal spot. If no such $v$
exists, then $\tau_{s_u} = \wmax{s_u}$ for every $u \ge 1$ and in particular $\tau_{s_1} = \wmax{s_1}$, then $s_1$ is
a maximal spot: there is no $s_0$, so the conditions at $s_{m-1}$ in
clause~(2) is vacuous.
\end{proof}
Next, we prove two propositions and one theorem about permissibility that were used in the proof of Theorem~\ref{thm:main_theorem}.
\begin{proposition} \label{prop:wiispermiss}
    Let $\sigma \in D_S \cap W(\tau)$, then $\bigl(w_1(\sigma), \dots,w_n(\sigma)\bigr)$ is permissible, thus 
    \begin{align}
         \bigcup_{\substack{\sigma\in S_n\\ \operatorname{Des}(\sigma)=S}}
\bigl(w_1(\sigma),\ldots,w_n(\sigma)\bigr)
\subseteq
\bigcup_{\tau\text{ permissible for }S}
(\tau_1,\ldots,\tau_n).
    \end{align}
\end{proposition}
\begin{proof}
    We have $\wmin{i} \leq w_i(\sigma) \leq \wmax{i}$. The lower bound $\wmin{i} \leq w_i(\sigma)$ follows from the fact that, in order to respect the descent set, $\sigma_i$ must be smaller than every entry after it in its own run. On the other hand, $w_i(\sigma) \leq \wmax{i}$ because the largest possible value of $w_i(\sigma)$ occurs when $\sigma_i$ is simultaneously larger than every entry in the next run.
    \\
    The second condition, namely $w_j(\sigma) \leq w_{j+1}(\sigma)+1$ when $j$ and $j+1$ lie in the same run, follows from the relative positions of $\sigma_j$ and $\sigma_{j+1}$. If $w_j(\sigma)=w_{j+1}(\sigma)+1$, then $\sigma_j$ and $\sigma_{j+1}$, with $\sigma_j<\sigma_{j+1}$, are larger than the same number of entries in the next run. If instead $w_j(\sigma)>w_{j+1}(\sigma)+1$, then $\sigma_j$ is larger than more entries in the next run than $\sigma_{j+1}$, which implies $\sigma_j>\sigma_{j+1}$. Since $j\notin S$, this contradicts the prescribed descent set. Therefore, we must have $w_j(\sigma)\leq w_{j+1}(\sigma)+1$.
\end{proof}
\begin{proposition}\label{prop:permissible}
    Let the descent set $S$ be given, and let $\tau \in \mathbb{N}^n$ be permissible, then $D_S \cap W(\tau) \ne \emptyset$, thus 
    \begin{align}
        \bigcup_{\substack{\sigma\in S_n\\ \operatorname{Des}(\sigma)=S}}
\bigl(w_1(\sigma),\ldots,w_n(\sigma)\bigr)
\supseteq
\bigcup_{\tau\text{ permissible for }S}
(\tau_1,\ldots,\tau_n).
    \end{align}
\end{proposition}
\begin{proof}
    Let $\tau$ be permissible. We construct a permutation $\sigma \in D_S \cap W(\tau)$ by inserting the entries $n,n-1,\dots,2,1$ successively.
    \\
    We begin by inserting $n$. Since $n$ can only be inserted at a maximal spot by Proposition~\ref{propmaxspot}, and a maximal spot always exists by Proposition~\ref{prop:maxspotexists}, it follows that $n$ can be inserted at a maximal spot.
    \\
    We then proceed by inserting the remaining entries in decreasing order,
$n-1,n-2,\ldots,1$. Suppose that we are inserting $y$, where $1\leq y<n$. At this stage, we modify the definition of a maximal spot to define a maximal spot while inserting $y$. A maximal spot encountered while inserting $y$ is called a $y$-maximal spot.
    \\
    Denote the rightmost available spot in each run $R_i$ by $r_i$. Write $I \subseteq \{1,\dots,d+1\}$ for the set of indices
$i$ such that $r_i$ is an available spot. Let $l = \min I$ and
$k = \max I$, so that $r_l$ and $r_k$ are the first and last available
spots. Note that $l$ and $k$ need not be $1$ and $d+1$. By the
definition of $w_i(\sigma)$, no two consecutive $r$'s are equal, so all
comparisons below are strict. A $y$-maximal spot is defined as follows:
    We call an available spot $r_i$ to be \emph{$y$-maximal} if at least one of
the following holds:
\begin{enumerate}
    \item $i = k$ and $r_i$ exceeds the preceding available spot;
    \item $i = l$ and $r_i$ exceeds the following available spot;
    \item one of:
    \begin{enumerate}
        \item $r_{i-1}$ and $r_{i+1}$ are both available, and
              $r_i > r_{i-1}$, $r_i > r_{i+1}$
              (a \emph{local maximum});
        \item neither $r_{i-1}$ nor $r_{i+1}$ is available;
        \item $r_{i-1}$ is not available and $r_{i+1} < r_i$;
        \item $r_{i+1}$ is not available and $r_{i-1} < r_i$.
    \end{enumerate}
\end{enumerate}
It remains to show that a $y$-maximal spot exists. The relative order of
consecutive available spots is determined by the corresponding $\tau$
values: $r_{i-1} + \tau_{r_{i-1}} < r_i$ gives $r_i > r_{i-1}$, and
$r_{i-1} + \tau_{r_{i-1}} > r_i$ gives $r_i < r_{i-1}$.\\
Suppose first that three consecutive available spots
$r_{i-1}, r_i, r_{i+1}$ exist. If some such $r_i$ satisfies
$r_i > r_{i-1}$ and $r_i > r_{i+1}$, it is a local maximum, hence
$y$-maximal by 3(a). If no local maximum exists, the available spots
form a monotone chain along $I$: either $r_l < \dots < r_k$, and $r_k$
is $y$-maximal by case~1, or $r_l > \dots > r_k$, and $r_l$ is
$y$-maximal by case~2.\\
If instead some neighbours are unavailable, the remaining
configurations are covered by 3(b)--(d). When neither $r_{i-1}$ nor
$r_{i+1}$ is available, $r_i$ is $y$-maximal by 3(b). When only
$r_{i+1}$ is available, either $r_{i+1} < r_i$ and $r_i$ is $y$-maximal
by 3(c), or $r_{i+1} > r_i$, in which case we follow the ascending chain
of available spots until it either turns down at some $r_p$, giving a
local maximum, or reaches $r_k$, which is $y$-maximal by case~1. The
situation with only $r_{i-1}$ available is symmetric, terminating either
at a local maximum or at $r_l$ by case~2.\\
Cases~1, 2 and 3(a)--(d) therefore exhaust all configurations of the
available spots, so at least one $y$-maximal spot exists.\\
At every stage, the condition $\tau_j \le \tau_{j+1} + 1$ for
consecutive elements $j, j+1$ of a run $R_i$ forces
$\sigma_j < \sigma_{j+1}$, so the resulting permutation respects the
prescribed descent set. The condition
$\wmin{i} \le \tau_i \le \wmax{i}$ for all $i$ ensures that the desired
$w_i$ values are realizable by inserting an element of $\{1,\dots,n\}$
while respecting $S$.\\
    Repeating this with the next largest value until all of
$\{1,\dots,n\}$ has been placed, we obtain a permutation with the
prescribed descent set and the prescribed values $w_i(\sigma)$.
\end{proof}
\begin{theorem} \label{thm:permissareallwsequence}
Let $n \in \N, S = \{s_1,\dots,s_d\} \subseteq \{1,\dots,n\}$ with $s_1 <\dots<s_d$. Then
    \begin{align}
\bigcup_{\substack{\sigma\in S_n\\ \operatorname{Des}(\sigma)=S}}
\bigl(w_1(\sigma),\ldots,w_n(\sigma)\bigr)
=
\bigcup_{\tau\text{ permissible for }S}
(\tau_1,\ldots,\tau_n).
    \end{align}
\end{theorem}
\begin{proof}
    This follows from Propositions~\ref{prop:wiispermiss} and~\ref{prop:permissible}.
\end{proof}
\noindent \textbf{Example 4.2.1} 
Construct a permutation with
\begin{align*}
    n = 10, S= \{2,4,7\}, \tau = (2,2,3,2,2,4,3,3,2,1)
\end{align*}
Note this $\tau$ is permissible for the descent set $S$. We use a vertical bar $|$ to denote a descent and a comma to denote a non-descent position in the illustration below.\\
First insert $10$. Condition 2 is satisfied in Definition 4.4, so 7 is a maximal spot.
\begin{align*}
    \bigl(\bigcirc, \bigcirc| \bigcirc, \bigcirc| \bigcirc, \bigcirc, 10| \bigcirc, \bigcirc, \bigcirc\bigr)
\end{align*}
Condition 3.(c) is satisfied in the definition of 9-maximal spot for 2, so 2 is a 9-maximal spot.
\begin{align*}
    \bigl(\bigcirc, 9| \bigcirc, \bigcirc| \bigcirc, \bigcirc 10| \bigcirc, \bigcirc, \bigcirc\bigr)
\end{align*}
Condition 3.(a) is satisfied in the definition of 8-maximal spot for 4, so 4 is a 8-maximal spot. 
\begin{align*}
    \bigl(\bigcirc, 9| \bigcirc, 8| \bigcirc, \bigcirc 10| \bigcirc, \bigcirc, \bigcirc\bigr)
\end{align*}
Condition 3.(c) is satisfied in the definition of 7-maximal spot for 1, so 1 is a 7-maximal spot.
\begin{align*}
    \bigl(7 \text{, } 9| \bigcirc, 8| \bigcirc, \bigcirc, 10| \bigcirc, \bigcirc, \bigcirc\bigr)
\end{align*}
Condition 3.(c) is satisfied in the definition of 6-maximal spot for 3, so 3 is a 6-maximal spot.
\begin{align*}
    \bigl(7 \text{, } 9| 6\text{, } 8| \bigcirc, \bigcirc, 10| \bigcirc, \bigcirc, \bigcirc\bigr)
\end{align*}
Condition 3.(c) is satisfied in the definition of 5-maximal spot for 6, so 6 is a 5-maximal spot.
\begin{align*}
    \bigl(7 \text{, } 9| 6\text{, } 8| \bigcirc,5, \text{ } 10| \bigcirc, \bigcirc, \bigcirc\bigr)
\end{align*}
Condition 1 is satisfied in the definition of 4-maximal spot for 10, so 10 is a 4-maximal spot.
\begin{align*}
    \bigl(7 \text{, } 9| 6\text{, } 8| \bigcirc,5 \text{, } 10| \bigcirc, \bigcirc, 4\bigr)
\end{align*}
Condition 1 is satisfied in the definition of 3-maximal spot for 9, so 9 is a 3-maximal spot.
\begin{align*}
    \bigl(7 \text{, } 9| 6\text{, } 8| \bigcirc,5, \text{ } 10| \bigcirc, 3 \text{, } 4\bigr)
\end{align*}
Condition 1 is satisfied in the definition of 2-maximal spot for 8, so 8 is a 2-maximal spot.
\begin{align*}
    \bigl(7 \text{, } 9| 6\text{, } 8| \bigcirc,5, \text{ } 10| 2, 3 \text{, } 4\bigr)
\end{align*}
Condition 3.(b) is satisfied in the definition of 1-maximal spot for 5, so 5 is a 1-maximal spot.\\
Now we have our permutation with the desired descent set and $w_i$ values
\begin{align*}
    \bigl(7 \text{, } 9| 6\text{, } 8| 1\text{, }5 \text{, } 10|2\text{, } 3 \text{, } 4\bigr)
\end{align*}
\begin{corollary}
    If $\tau$ is permissible for the descent set $S$, then $|D_S \bigcap W(\tau, U)| \ne 0$ for any $U \in 2^{\{1,\dots,s_d\}}$.
\end{corollary}
This follows immediately from the inclusion $W(\tau,\phi)\subseteq W(\tau,U)$ for every $U\ne\emptyset$. Since $\tau$ is permissible for $S$, we have $D_S\cap W(\tau,\phi)\ne\emptyset$, and hence $D_S\cap W(\tau,U)\ne\emptyset$.\\
We are now ready to state and prove the main counting theorem, which provides the formula for $|D_S\cap W(\tau,U)|$ that was used in the proof of Theorem~\ref{thm:main_theorem}.
\begin{theorem} \label{maincount}
    Let $n \in \N, S=\{s_1, s_2,\dots, s_d\} \subseteq \{1,\dots, b\}$ with $s_1 < \dots < s_d$ and $s_1 + s_2 +\dots + s_d = b$. Let $U \subseteq \{1,2,\dots, s_d\}$ and $\tau \in \mathbb{N}^n$ be permissible for $S$, then $|D_S \cap W(\tau,U)|$ is a polynomial in $n$.
\end{theorem}
\subsection{Proof of Theorem \ref{maincount}} \label{sec:mainproof}
We prove by using recurrence relations in $|D_S \cap W(\tau,U)|$. Before presenting the proof of the recurrence relation, we illustrate its use by deriving a polynomial formula for $|D_S \cap W(\tau,U)|$ with $S=\{1,3,5\}$, $U=\{5\}$, and $\tau=(1,2,2,1,3)\sqcup(n-5)!$. To facilitate this computation, we first establish the following three lemmas, a corollary and a new notation.
\begin{lemma}[Telescoping]\label{lem:telescope}
Let $N \in \mathbb{Z}$ and let $(F_n)_{n \ge N-1}$ be a sequence of rational
numbers satisfying
\begin{enumerate}
  \item $F_{N-1} = 0$, and
  \item $F_n = F_{n-1} + Q(n)$ for all $n \ge N$,
\end{enumerate}
where $Q$ is any function on the integers $\ge N$. Then
\[
  F_n \;=\; \sum_{i=N}^{n} Q(i) \qquad \text{for all } n \ge N-1,
\]
with the convention that the empty sum (the case $n = N-1$) is $0$.
\end{lemma}

\begin{proof}
Fix $n \ge N-1$. If $n = N-1$ both sides are $0$ by condition 1. Otherwise
$n \ge N$, and condition 2 applies at each index $i$ with
$N \le i \le n$, giving $F_i - F_{i-1} = Q(i)$. Summing over this range,
\[
  \sum_{i=N}^{n}\bigl(F_i - F_{i-1}\bigr) \;=\; \sum_{i=N}^{n} Q(i).
\]
The left-hand side telescopes: each $F_i$ with $N \le i \le n-1$ occurs
once positively and once negatively, leaving $F_n - F_{N-1} =F_n = \sum_{i=N}^{n} Q(i)$ by condition 1.
\end{proof}
\begin{lemma}[Power sums]\label{lem:powersum}
For all integers $p \ge 0$ and $n \ge 0$,
\[
  \sum_{i=1}^{n} i^{p} \;=\; \sum_{j=1}^{p+1} \frac{1}{j}\, Stir(p+1,j)\, (n)_j ,
\]
where $Stir(a,b)$ is the Stirling number of the second kind and
$(n)_j = n(n-1)\cdots(n-j+1)$ is the falling factorial. A proof can be found in \cite{GKP1994}.
\end{lemma}
\begin{corollary}\label{cor:polysum}
Let $Q \in \mathbb{Q}[x]$ have degree $d$ and let $m \in \mathbb{Z}$ be
fixed. Then there is a unique $R \in \mathbb{Q}[x]$ of degree $d+1$ with
\[
  \sum_{i=m}^{n} Q(i) \;=\; R(n) \qquad \text{for all } n \ge m-1 .
\]
\end{corollary}

\begin{proof}
Each $(n)_j$ is a product of $j$ linear factors, hence a polynomial in $n$
of degree $j$; the sum in Lemma~\ref{lem:powersum} has $p+1$ terms, so
$\sum_{i=1}^{n} i^{p}$ is a polynomial in $n$, of degree $p+1$ since only
the $j = p+1$ term contributes that degree. Writing
$Q(x) = \sum_{j=0}^{d} a_j x^{j}$ and summing termwise,
$\sum_{i=1}^{n} Q(i)$ is therefore a polynomial in $n$ of degree $d+1$. Finally
\[
  \sum_{i=m}^{n} Q(i) \;=\; \sum_{i=1}^{n} Q(i) \;-\; \sum_{i=1}^{m-1} Q(i),
\]
and the second sum is a constant independent of $n$, so the difference is
again a polynomial of degree $d+1$. Uniqueness holds because two
polynomials agreeing at infinitely many integers are equal.
\end{proof} 
\begin{notation}\label{not:bar}
Let $k \le n$, let $\tau = (\tau_1,\dots,\tau_k) \sqcup (n-k)!$, and let
$S, U \subseteq [k]$. We abbreviate the set $D_S \cap W(\tau,U)$ by a
single string: write the entries $\tau_1,\dots,\tau_k$ in order, followed
by $(n-k)!$; separate the entries in positions $i$ and $i+1$ by a bar
$\mid$ if $i \in S$ and by a comma otherwise; and replace $\tau_i$ by
$\ge \tau_i$ if $i \in U$. For example, taking $S = \{1,3,5\}$,
$U = \{5\}$ and $\tau = (1,2,2,1,3) \sqcup (n-5)!$, our notation gives
\[
  D_{\{1,3,5\}} \cap W\bigl((1,2,2,1,3) \sqcup (n-5)!,\ \{5\}\bigr)
  \;=\; \bigl(1 \mid 2, 2 \mid 1, \ge 3 \mid (n-5)!\bigr),
\]
and we write $\#\bigl(1 \mid 2, 2 \mid 1, \ge 3 \mid (n-5)!\bigr)$ for its
cardinality.

Thus a bar records a descent, a comma records a non-descent, the symbol
$\ge$ records the truncated term, and $(n-k)!$ records the last run of $\tau$ which by permissibility has to be $(n-k, n-k-1, \dots, 2,1)$. The string determines $S$, $U$ and
$\tau$ uniquely, so the notation loses no information. We reserve
$\#$ rather than $|\cdot|$ for cardinality throughout, to avoid confusion
with the descent bars.
\end{notation}
\begin{lemma}\label{lem:base}
Let $\sigma$ range over $S_m$ and suppose $S = \{1\}$. Then
$\#\bigl(\ge \tau_1 \dsc (m-1)!\bigr) = m - \tau_1$ and
$\#\bigl(\tau_1 \dsc (m-1)!\bigr) = 1$.
\end{lemma}
\begin{proof}
    For $U = \{1\},$ allowing $\omega_1(\sigma)$ be any number greater than $\tau_1$ means that $\sigma_1 > \sigma_2$, $\sigma$ is increasing starting from $\sigma_2$, and that the first $\tau_1$ entries in the second run are smaller than $\sigma_1$. To construct such a permutation, we need $\sigma_1 > \tau_1$, and the rest of the permutation is just the remaining entries put in increasing order. For example $(\ge 2|3,2,1)$ consists of permutations $(3,1,2,4)$ and $(4,1,2,3)$. \\
    For $U = \emptyset$, the only permutation that satisfies the condition is $(2,1,3,4,...,m-1, m)$.
\end{proof}
Now we calculate $\#\bigl(1 \mid 2, 2 \mid 1, \ge 3 \mid (n-5)!\bigr)$. We first introduce the tree structure.
\paragraph{Branching.}
Each child of a node is obtained by deleting the largest entry of the
parent and recording how $S$, $\tau$ and $U$ change. The cases are
indexed by the position of that entry (which are maximal spots), shown on the edges.
Definitions~\ref{phimap} and~\ref{psimap} will give the maps describing these changes ($\phi_t$ for $S$, $\psi_t^S$ for $\tau$ and
$U$), and Theorem~\ref{verphipsi} verifies them; the reader may
also check any individual edge directly.

\paragraph{Terminal nodes of the first kind.}
We stop at a node whose descent set is $\{1\}$, since
Lemma~\ref{lem:base} evaluates it outright.

\paragraph{Terminal nodes of the second kind.}
We stop at the child obtained by deleting $n-i$ from the \emph{last} position at generation $i$. This deletion leaves $\tau_1, \tau_2, \dots,\tau_{s_d}$, $S$, and $U$ unchanged, while shortening the last run by one. Let $F_m$ denote the value of a node corresponding to permutations of length $m$ with fixed $S$, $U$, and $\tau_1,\tau_2,\dots,\tau_{s_d}$. Then the child obtained by this deletion has value $F_{m-1}$. This gives the recurrence step, which we carry out while constructing the tree.
\begin{center}
\begin{forest}
  for tree={
    grow=south,
    s sep=50pt,
    l sep=30pt,
  }
  [{$\#\bigl(1 \dsc 2 \asc 2 \dsc 1 \asc \ge 3 \dsc (n-5)!\bigr)$}
    [{$\#\bigl(1 \dsc 1 \dsc 1 \asc \ge 3 \dsc (n-5)!\bigr)$},
      edge label={node[midway,fill=white,font=\scriptsize]{$\sigma_3 = n$}}
      [{$\#\bigl(1 \dsc 1 \dsc 1 \asc \ge 3 \dsc (n-6)!\bigr)$},
        edge label={node[midway,fill=white,font=\scriptsize]{$\sigma_{n-1} = n-1$}}]
      [{$\#\bigl(1 \dsc \ge 1 \dsc (n-4)!\bigr)$},
        edge label={node[midway,fill=white,font=\scriptsize]{$\sigma_4 = n-1$}}]
      [{$\#\bigl(1 \dsc 1 \asc \ge 3 \dsc (n-5)!\bigr)$},
        edge label={node[midway,fill=white,font=\scriptsize]{$\sigma_1 = n-1$}}
        [{$\#\bigl(1 \dsc 1 \asc \ge 3 \dsc (n-6)!\bigr)$},
          edge label={node[midway,fill=white,font=\scriptsize]{$\sigma_{n-2} = n-2$}}]
        [{$\#\bigl(\ge 1 \dsc (n-4)!\bigr)$},
          edge label={node[midway,fill=white,font=\scriptsize]{$\sigma_3 = n-2$}}]
      ]
    ]
    [{$\#\bigl(1 \dsc 2 \asc 2 \dsc 1 \asc \ge 3 \dsc (n-6)!\bigr)$},
      edge label={node[midway,fill=white,font=\scriptsize]{$\sigma_n = n$}}]
  ]
\end{forest}
\end{center}
We begin the calculation from the bottom. For $\sigma\in\bigl(\ge1\dsc(n-4)!\bigr)$, we have $S=\{1\}=U$, so this is a base case. Moreover, $\sigma \in S_{n-3}$. Thus, by Lemma~\ref{lem:base}, $\#\bigl(\ge 1 \dsc (n-4)!\bigr)= (n-3)-1 = (n-2) -2$. If we denote $F_{n-2} = \#\bigl(1 \dsc 1 \asc \ge 3 \dsc (n-5)!\bigr)$, then we have the recurrence step 
\begin{equation}
    F_{n-2} = F_{n-3} + (n-2) - 2
\end{equation}
Letting $m = n-2$, this is equivalent to
\begin{align} \label{frecur}
    F_m = F_{m-1} + m - 2 
\end{align}
Since $F_m$ counts $\#\bigl(1 \dsc 1 \asc \ge 3 \dsc (m-3)!\bigr)$, we know that at $F_5 = 0$ (it is impossible to have $\tau_3 \ge 3$ while the next run has only length 2) and the recurrence (\ref{frecur}) starts at $F_6$. By Lemma \ref{lem:telescope} and \ref{lem:powersum}
\begin{align}
    F_m = \sum_{i=6}^m (i-2) &= \frac{1}{2}m^2 - \frac{3}{2}m - 5\\
\end{align}
Converting back to $n$ we get:
\begin{align}
    F_{n-2} &= \frac{1}{2}(n-2)^2 - \frac{3}{2}(n-2) - 5\\
    &= \frac{1}{2}n^2 - \frac{7}{2}n
\end{align}
We should have two branches after
$\#\bigl((1\mid\geq 1\mid(n-4)!)\bigr)$, with the non-terminal node of
the second kind reaching the base case $S=\{1\}$. However, this
node is straightforward to count directly: it is precisely the number
of permutations of $n-2$ with descent set ${1,2}$. We therefore count
it directly here.\\
To construct such a permutation, it suffices to choose two distinct
numbers from ${2,3,\ldots,n-2}$ to be $\sigma_1$ and $\sigma_2$,
with the smaller of the two assigned to $\sigma_2$, and then append the
remaining numbers in increasing order after $\sigma_2$. Hence,
\(
\#\bigl((1\mid\geq 1\mid(n-4)!)\bigr)=
\binom{(n-2)-1}{2}.
\)
We would then apply the recurrence to
$\#\bigl(1\mid1\mid1,\geq3\mid(n-4)!\bigr)$ in a similar fashion. We
omit the remaining calculations here and instead provide the final
tree of polynomials.

\begin{center}
\begin{forest}
  for tree={
    grow=south,
    s sep=50pt,
    l sep=30pt,
  }
  [{$\frac{1}{12}n^4 - \frac{5}{6}n^3 - \frac{1}{12}n^2 + \frac{89}{6}n - 14$}
    [{$\frac{1}{3}n^3 - 3n^2 - \frac{8}{3}n + 14$},
      edge label={node[midway,fill=white,font=\scriptsize]{$\sigma_3 = n$}}
      [{recurrence step},
        edge label={node[midway,fill=white,font=\scriptsize]{$\sigma_{n-1} = n-1$}}]
      [{$\binom{(n-2)-1}{2}$},
        edge label={node[midway,fill=white,font=\scriptsize]{$\sigma_4 = n-1$}}]
      [{$\frac{1}{2}n^2 - \frac{7}{2}n$},
        edge label={node[midway,fill=white,font=\scriptsize]{$\sigma_1 = n-1$}}
        [{recurrence step},
          edge label={node[midway,fill=white,font=\scriptsize]{$\sigma_{n-2} = n-2$}}]
        [{$(n-3)-1$},
          edge label={node[midway,fill=white,font=\scriptsize]{$\sigma_3 = n-2$}}]
      ]
    ]
    [{recurrence step},
      edge label={node[midway,fill=white,font=\scriptsize]{$\sigma_n = n$}}]
  ]
\end{forest}
\end{center}
Having presented the example, we will now establish that the results observed hold true in general. We start by giving details to how $S, \tau$ and $U$ change by removing the biggest entry. First, we study how the descent set $S$ change.
\begin{definition} \label{phimap}
    For $ \sigma_t = s_y$ a maximal spot, define $\phi_t: 2^{\{1,...,b\}} \to 2^{\{1,...,b\}}$ as follow:
\begin{enumerate}
    \item $t = n$. Then $\phi_t(S) = S$
    \item $1 < t < n$. 
    \begin{enumerate}
        \item If $\tau_{t-1} = 1,$ then $s_i \mapsto s_i - 1$ for $i > y$, $s_i \mapsto s_i$ for $i<y$, and then delete $m = s_y$. 
        \item If $\tau_{t-1} \ge 2,$ then $s_i \mapsto s_i$ for $1 \le i < y$, $s_i \mapsto s_i - 1$ for $i \ge y$.
    \end{enumerate}
    \item $t = 1$, $s_i \mapsto s_i - 1$ and then delete $s_1$.
\end{enumerate}
\end{definition}
Now we define $\psi$, which tracks the changes in $\tau$ and $U$ under the removal of $\sigma_t=n$. Observe that if $\sigma_t=n$ with $t\in S$, then we can compare $\sigma_{t-1}$ and $\sigma_{t+1}$. If $\tau_{t-1}=1$ and $t-1\notin U$, then $\sigma_{t-1} < \sigma_{t+1}$, whereas if $\tau_{t-1}\ge2$, then $\sigma_{t-1} > \sigma_{t+1}$.
\begin{definition} \label{psimap}
    Let $S$ be given, $U \subseteq \{1,...,s_d\}$, and $\tau \in \N^n$ permissible. For $\sigma_t = s_y$ a maximal spot, define $\psi_t^S:  \N^n \times 2^{\{1,2,...,n\}} \to \N^{n-1} \times 2^{\{1,2,...,n-1\}}$ as follows:
\begin{enumerate}
    \item $t = n$\\
    First delete $w_n$, and for $i \in \{s_d, s_d+1,..., n-1\}$, $w_i \mapsto w_i - 1$. $w_i$ stays the same for all other $i$. $U$ stays the same.
    \item $1 < t < n$. 
    \begin{enumerate}
        \item $\tau_{t-1} = 1$ and $t-1 \notin U$.\\
\emph{Effects on $\tau$}: Recall that the $i$-th run $R_i$ denotes the $i$-th increasing sequence, then $t = s_y$ is in $R_y$.
\begin{itemize}
  \item If $y < d$, then for each $i \in R_y$,
        \[ \tau_i \;\longmapsto\; \tau_i - 1 + (s_{y+1} - s_y). \]
  \item If $y = d$, then $\tau_i \mapsto n - i$ for each $i \in R_d$.
\end{itemize}
All $\tau_i$ with $i \notin R_y$ are unchanged.

\emph{Effects on $U$}: Suppose $t \in U$. Then $U$ changes in two ways:
\begin{itemize}
  \item every $u \in U$ with $u > t$ is reindexed as $u - 1$;
  \item every $i \in R_{y-1}$ with $\tau_i + i \ge t - 1$ is added.
\end{itemize}
Writing $U'$ for the result,
\[
  U' \;=\; \{u : u \in U,\ u < t\}
        \;\cup\; \{u - 1 : u \in U,\ u > t\}
        \;\cup\; \{i \in R_{y-1} : \tau_i + i \ge t-1\}.
\]
        \item $\tau_{t-1} \ge 2$\\
        First, delete $\tau_t$. For each $i$ in the same run as $m$, i.e.,
$i\in\{s_{y-1},s_{y-1}+1,\ldots,s_y\}$, replace $\tau_i$ by
\begin{align}
\begin{cases}
\tau_i-1, & \text{if $\tau_i>1$},\\
1, & \text{if $\tau_i=1$}.
\end{cases}
\end{align}
For the index set $U$, retain each $u\in U$ with $u<t$ unchanged, while
replacing each $u\in U$ with $u-1$ for $u>t$.
    \end{enumerate}
   \item $m=1$.\\
Delete $\tau_1$. If $1\in U$, delete $1$ from $U$ and replace each
remaining $u\in U$ by $u-1$.
\end{enumerate}
If $\tau_{t-1} = 1$ and $t - 1 \in U$, then we need to treat $\psi\big(W(\tau, U)\big)$ as $\psi\big(W(\tau', U)\big) \bigcup \psi\big(W(\tau, U')\big)$ where $\tau'_i = \tau_i$ for $i \ne t - 1$ and $\tau_{t-1} = 2$, and $U' = U - \{t-1\}$, since our map have different treatment for $\tau_{t-1}=1$ and $\tau_{t-1} \ge 2$.
\end{definition}
\begin{definition}\label{def:Z}
Denote $[x] := \{1, 2,\dots, x-1,x\}.$ Let $\sigma \in S_n$ and let $t \in [n]$ be the position of the
largest entry of $\sigma$. Combining Definitions~\ref{phimap}
and~\ref{psimap}, define
\[
  Z_t : 2^{S_n} \to 2^{S_{n-1}}
\]
by
\[
  Z_t(D_S \cap W(\tau,U))
  \;=\;
  D_{\phi_t(S)} \cap W\bigl(\psi_t^{S}(\tau, U)\bigr).
\]

\end{definition}
We now work through two examples in detail. These are the first two
generations of the tree above, with $n = 10$: we delete the entries $10$
and $9$ in turn and record precisely how the descent set $S$ and the
values $w_i$ change.
\begin{example}
    Let $n = 10, U = \{5\}$, $S = \{1,3,5\}, \tau = ( 1,2,2,1,3,5,4,3,2,1)$. In other words we are looking for permutations in $(1|2,2|1,3|5,4,3,2,1)$.\\
First we need to find the maximal spots. $\tau_3 = 2 = \wmax{3}$ so 3 is a maximal spot. Since $\tau_5 < 5$, we have that $10$ is also a maximal spot. Therefore we have in total of two maximal spots.
\begin{enumerate}
    \item $t = 3$\\
    Since $\tau_2 = 2$ we know that $\sigma_2 > \sigma_4$, so following Definition \ref{phimap}.2.(b), we have that $\phi_{3}(S) = \{1,2,4\}$.\\
    By following Definition \ref{psimap}.2(b), we have 
    \begin{align}
        Z_3\big((1|2,2|1,3|5,4,3,2,1)\big) = (1|1|1,\ge3|5,4,3,2,1)
    \end{align}
    \item $t = 10$\\
    By following Definition \ref{phimap}.1, we have that $\phi_{10}(S) = \{1,3,5\}$.\\
    By following Definition \ref{psimap}.1, we have that $\psi_{10}\big(W(\tau,U)\big) = W\big((1,2,2,1,3,4,3,2,1), \{5\}\big)$.
    \begin{align}
        Z_{10}\big((1|2,2|1,3|5,4,3,2,1)\big) = (1|2,2|1,\ge 3|4,3,2,1)
    \end{align}
\end{enumerate}
\end{example}
\begin{example}
    Let $n=9, S=\{1,2,4\}, U = \{5\}, \tau = (1,1,1,3,5,4,3,2,1)$. In other words we are looking for permutations of 9 in $(1|1|1, \ge3 |5,4,3,2,1)$. \\
First we need to find the maximal spots. 1 is a maximal spot, since $\tau_1 = \wmax{1} = 1$. 4 is also a maximal spot, and 9 is a maximal spot when $\tau_4 <5$. 
\begin{enumerate}
    \item $t = 1$\\
    By following Definition \ref{phimap}.3 and \ref{psimap}.3, we get that $\phi_9(S) = \{1,3\}$ and 
    \begin{align}
        Z_1\big((1|1|1, \ge3 |5,4,3,2,1)\big) = (1|1, \ge3|5,4,3,2,1)
    \end{align}
    \item $t = 4$\\
    Since $4 \in U$, we can have $\tau_4 = \wmax{4}(\sigma) = 5$, and $\sigma_4$ is a maximal spot. Note that since $\tau_2 = 1$ we have that $\sigma_2< \sigma_4$, so if we follow Definition \ref{phimap}.2(a) and \ref{psimap}.2(a), we have that $\phi_4(S) = \{1,2\},$ and
    \begin{align}
        Z_4\big((1|1|1, \ge3 |5,4,3,2,1)\big) = (1|\ge1|6,5,4,3,2,1)
    \end{align}
    \item $t = 9$\\
    Since $4 \in U$ implies $\tau_4 \ge 3$, in this case $\tau_4 =3$ or 4, so $t = 9$ satisfies Definition~\ref{def:maxspot}.1 and is a maximal spot. By following Definition \ref{phimap}.3 and \ref{psimap}.3, we have that $\phi_9(S) = \{1,2,4\}$, and
    \begin{align}
        Z_9\big((1|1|1, \ge3 |5,4,3,2,1)\big) = (1|1|1, \ge3 |4,3,2,1)
    \end{align}
\end{enumerate}
\end{example}
\begin{theorem} \label{verphipsi}
    Let $S$ be given, $\tau$ be permissible and $U \subseteq \{1,2,...,s_d\}$. Then 
    \begin{align}
    \bigg|D_S \bigcap W(\tau, U)\bigg| = \sum_{m \text{ a maximal spot}} \bigg|D_{\phi_t(S)} \bigcap \psi_t^S \bigl(W(\tau,U)\bigl)\bigg| 
\end{align}
\end{theorem}
\begin{proof}
    Let $\sigma \in D_S \bigcap W(\tau, U).$ First, we will prove that if $\sigma_t = n$ where $t$ is a maximal spot, then removing $n$ yields precisely a permutation $\sigma' \in D_{\phi_t(S)} \bigcap \psi_t^S \big(W(\tau;U)\big)$. \\
    We first show that $\sigma' \in D_{\phi_t(S)}$. In case 1, deleting the
last entry does not alter the descent set. In case 2.(a), the condition
$\sigma_{t-1}<\sigma_{t+1}$ implies that deleting
$\sigma_t=n$ merges the two runs $R_y$ and $R_{y+1}$. Thus, the descent
at $s_y$ is removed from the descent set, and the indices of all
subsequent descents are shifted accordingly. In case 2.(b), deleting
$n=\sigma_t$ does not merge the two runs, so the only change is the
corresponding shift in the indices of the subsequent descents. Finally,
in case 3, deleting $\sigma_1=n$ does not affect any subsequent descent,
and hence only the indices of the descents are shifted.\\
    Now we prove that $\sigma' \in \psi_t^S \big(W(\tau,U)\big).$ We divide into 3 cases as before.
\begin{enumerate}
    \item If $t=n$, then no elements in the previous run is affected by removing $n$, but every element in the same run as $n$ is less than $n$, so removing $n$ would have $w_i$ decreases 1 as defined. No other elements would have $w_i$ value changed by removing $n$.
    \item $1 < t < n$. Suppose $t = s_y$
    \begin{enumerate}
        \item $\sigma_{t-1} < \sigma_{t+1}$\\
        For elements $\sigma_i$ in the same run $R_y$ as $\sigma_t=n$, removing $n$ removes an element larger than $\sigma_i$, so we obtain $w_i-1$. In the new permutation, since $\sigma_{t+1}>\sigma_{t-1}$, the next run becomes the same run as $\sigma_j$. Moreover, by the descent set, every element in the next run is larger than $\sigma_j$. Thus, we must add the length of the next run, giving
\(
w_i-1+(s_{y+1}-s_y).
\)
When $y=d$, we modify the map accordingly, using the placeholder term $\sigma_{n+1}=0$.\\
        For elements $\sigma_i$ in the previous run $R_{y-1}$ of $n$, suppose $\sigma_i \in R_{y-1}$ satisfies
$i + w_i(\sigma) \ge t-1$, so that $\sigma_i > \sigma_{t-1}$. Removing
$n$ merges $R_y$ and $R_{y+1}$, so the run following $\sigma_i$ may now
contain additional elements smaller than $\sigma_i$.

In fact every value in the range is attained: for each $z$ with
$0 \le z \le s_{y+1}-s_y$, if exactly $z$ of
$\sigma_{s_y},\dots,\sigma_{s_{y+1}}$ are smaller than $\sigma_i$, then
reinserting $n$ at position $t$ yields a permutation in
$D_S \cap W(\tau,U)$ with $w_i(\sigma) = z$. All such elements of the
preceding run therefore lie in $U$, which is what the definition of
$\psi_t^S$ records in this case.\\
        If instead $i + w_i(\sigma) < t-1$, then every element of the run
following $n$ is greater than $\sigma_i$. The merged run is longer, but
the elements it gains all exceed $\sigma_i$, so the number of elements
in the following run that are smaller than $\sigma_i$ is unchanged.
This is what the definition of $\psi_t^S$ records in this case.
        \item $\sigma_{t-1} > \sigma_{t+1}$\\
        For elements in the run preceding that of $n$, the $w_i$ values are
unchanged: removing $n$ does not merge the two following runs, and does
not increase the number of elements in $n$'s run that are smaller than
any given element of the preceding run. The only values affected are
those of elements in the same run as $n$, for which $w_i \mapsto w_i-1$.
This is what the definition of $\psi_t^S$ records in this case.
    \end{enumerate}
    \item $t = 1$\\
    Suppose $m = 1$, so that $\sigma_1 = n$. No element after $n$ has its
$w_i$ value affected by the removal, so $w_1$ is deleted and the
remaining entries are reindexed, with $U$ relabelled accordingly as
$\{u - 1 : u \in U\}$.
\end{enumerate}
Therefore we have proven that $\sigma' \in \psi_t^S \bigl(W(\tau;U)\bigl)$. Since we have proven that $n$ appear and only appear at the maximal spots, we get the result by summing over all the maximal spots.\\
\end{proof}
Now we are ready to give a proof of Theorem \ref{maincount}.
\begin{proof}
    First we formalize the action of consecutively applying $Z_t$ multiple times. 
    Let $Z^j_{(t_1,...,t_j)}( S, \tau, U)$ denote 
$Z_{t_j} \circ Z_{t_{j-1}}...\circ Z_{t_1}(S, \tau, U)$, where $t_i$ is a maximal spot of $Z^{i-1}_{(t_1,...,t_{i-1})}(S, \tau, U)$ for all $i$. For example, in the tree we have shown in the beginning of this section, we have that $Z_{(3,4)}^2\Big(\big((1|2,2|1,\ge3|(n-5)!)\big)\Big) = \big((1|\ge1|(n-4)!\big).$\\
First, prove that $|Z^1_{t_1}(S, \tau, U)|$ where $t_1 \ne n$ is a polynomial in $n$ ($t_1 \ne n$ ensures that this branch is not the aforementioned terminal node of the second kind). 
By nature of the map $\phi$, every time $\phi$ is applied the last descent position $s_d$ decreases at least by one, so we know that by at most $s_d - 1$ times of applying $\phi$, we can get to the descent set $S = \{1\}. $ Let the number of times needed to be applied to $S$ to get to $S = \{1\}$ be $j$. We will prove the aforementioned statement by reverse induction on the number of times $Z$ map is applied.\\
First, prove that for the base case $|Z^j_{(t_1,...,t_j)}(S, \tau, U)|$ where $t_1,\dots, t_j$ is the sequence of maximal spots that lead to the base case $S = \{1\}$. Then $|Z^j_{(t_1,...,t_j)}(S, \tau, U)|$ is a polynomial in $n$. Indeed we know that the final descent set is $\{1\}$, and by Lemma \ref{lem:base} this is automatically true.\\
For the induction step,  assume that $|Z^i_{(t_1,...,t_i)}(S, \tau, U)|$, where $t_1 \ne n, t_2 \ne n - 1,..., t_i \ne n - i + 1$,  is a polynomial in $n$. The inequalities ensure that none of the nodes passed were terminal node of the second kind. Prove that $|Z^{i-1}_{(t_1,...,t_{i-1})}(S, \tau, U)|$, where $t_1 \ne n, t_2 \ne n - 1,..., t_{i-1} \ne n - i$ is a polynomial in $n$.\\
By Theorem~\ref{verphipsi},
\[
  \bigl|Z^{i-1}_{(t_1,\dots,t_{i-1})}(S,\tau,U)\bigr|
  = \sum_{t_i^k} \Bigl|Z_{t_i^k}\bigl(Z^{i-1}_{(t_1,\dots,t_{i-1})}(S,\tau,U)\bigr)\Bigr|,
\]
the sum being over the maximal spots $t_i^k$.\\
Suppose first that $t_i^k \ne n-i+1$ (ensuring $t_i^k$ is not the right most position in the permutation thus avoiding producing a terminal node of the second kind), then each
$\bigl|Z_{t_i^k}\bigl(Z^{i-1}_{(t_1,\dots,t_{i-1})}(S,\tau,U)\bigr)\bigr|
= \bigl|Z^{i}_{(t_1,\dots,t_{i-1},t_i^k)}(S,\tau,U)\bigr|$
falls under the induction hypothesis, so each such term is a polynomial
in $n$. Their sum, which we denote $P(n)$, is again a polynomial: the
number of summands is the number of maximal spots, which is finite and
depends only on $s_d$, hence only on $b$, not on $n$. Rewriting $P(n)$
in terms of $n-i+1$, put $P(n) = Q(n-i+1)$.\\
Now suppose $t_i^k = n-i+1$. By parts~1 of Definitions~\ref{phimap}
and~\ref{psimap}, applying $Z_{n-i+1}$ does not change $S, U$ and the only change is the length of the last run of $\tau$ to being one shorter. For instance, if the
last run of $Z^{i-1}_{(t_1,\dots,t_{i-1})}(S,\tau,U)$ contributes
$(n-p)!$, then the corresponding contribution after applying
$Z_{n-i+1}$ is $(n-p-1)!$.
This gives the recurrence:
\begin{align} \label{equabiglong}
    \big|Z^{i-1}_{(t_1,...,t_{i-1})}(S, \tau, U)\big| = Q(n - i +1) + \big|Z^{i}_{(t_1,...,t_{i-1}, t_i)}(S, \tau, U)\big|
\end{align}
If we denote $\big|Z^{i-1}_{(t_1,...,t_{i-1})}(S, \tau, U)\big| $ with total length of its w-sequence being $n-i+1$ as $F_{n-i+1}$, and $\big|Z^{i}_{(t_1,...,t_{i-1}, t_i)}(S, \tau, U)\big|$ with total length of its w-sequence being $n-i$ as $F_{n-i}$, then (\ref{equabiglong}) can be rewritten as 
\begin{align}
    F_{n-i+1} = F_{n-i} + Q(n-i+1)
\end{align}
Letting $m = n-i+1$, this is equivalent to 
\begin{align}
    F_m = F_{m-1} + Q(m)
\end{align}
Let $\alpha$ denote the last descent at this node, then we know the recurrence starts at $F_{\alpha + \tau_\alpha}$ and set $F_{\alpha + \tau_\alpha-1} = 0$. By Lemma~\ref{lem:powersum}, ~\ref{lem:telescope} and Corollary~\ref{cor:polysum}
\begin{align}
    F_m = \sum_{i=\alpha+ \tau_\alpha}^m Q(i)
\end{align}
Substituting back $m = n-i+1$ we get
\begin{align}
    F_{n-i+1} = \big|Z^{i-1}_{(t_1,...,t_{i-1})}(S, \tau, U)\big| = \sum_{l = \alpha + \tau_\alpha}^{n - i+1} Q(l)
\end{align}
Note that this gives a positive value precisely when
$n-i+1\geq \alpha+\tau_\alpha$. When $n-i+1=\alpha+\tau_\alpha$,
the value is $0$, which is also the correct value in this case.
\\
Recall that the maximum number of times the map $Z$ needs to be applied to get to the base case was $s_d-1$ thus $b-1$, so by Corollory~\ref{cor:polysum} and reverse induction, we have proven that $Z^1_{t_1}(S, \tau, U)$ where $t_1 \ne n$ is a polynomial in $n$. Using similar reasoning, we can prove that $\big|D_S \bigcap W(\tau,U)\big| = Z^0(S, \tau, U)$ is a polynomial in $n$. We know that 
\begin{align*}
    |D_S \bigcap W(\tau, U)|  = \sum_{\text{maximal spots } t_1^k} |Z_{t_1^k}(S, \tau, U)|
\end{align*}
For $t_1^k \ne n,$ we have that 
$\sum_{t_1^k \ne n} |Z_{t_1^k}(S, \tau, U)|$
is a polynomial in $n - 1$ denoted as $P(n-1)$, and after rewriting it into a polynomial of $n$ we get a polynomial $Q(n)$. For $t_1 = n$, we get that $Z_n(S, \tau, U)$ does not change $S$ and $U$ and only change the last run of $\tau$ to be one shorter, so we get 
\begin{align*}
    |D_S \bigcap W(\tau, U)| = Z^1_n(S, \tau, U) + Q(n)
\end{align*}
By similar reasoning as before, we can find that the function that satisfies this recurrence is a polynomial in $n$, namely
\begin{align*}
    |D_S \bigcap W(\tau, U)| = \sum_{l = s_d + \tau_{s_d}}^n Q(l)
\end{align*}
Now we have proved the theorem.\\
    \renewcommand{\qedsymbol}{$\blacksquare$}
\end{proof}
\begin{remark} \label{permissvsU}
When applying the formula for $D_S\cap W(\tau,U)$, it may happen that
the $\tau$ obtained by requiring 
$\tau_l\leq k$ for $l \not\in U$ is not
permissible because $\tau_l < \wmin{l}$. In this case, we dissect the configuration according to
the possible values of one of the entries until the resulting
configurations have permissible underlying $\tau$.
For example for $k = 0$, we may encounter calculating $\#(\geq 1,\geq 1,\geq 1\mid 3,2,1)$ while $(1,1,1,3,2,1)$ is not permissible, but $(\geq 1,\geq 1,\geq 1\mid 3,2,1)$ is not empty. In this case we dissect the union into permissible and non-permissible parts
\begin{align}
\#(\geq 1,\geq 1,\geq 1\mid 3,2,1)
&=
\#(1,\geq 1,\geq 1\mid 3,2,1)
+
\#(\geq 2,\geq 1,\geq 1\mid 3,2,1)\\
&=
\#(\geq 2,\geq 1,\geq 1\mid 3,2,1),
\end{align}
like so, and in this case
$(1,1,1,3,2,1)$ is not permissible. The remaining term 
$(\ge2,1,1|3,2,1)$ however is permissible, and the formula for
$D_S\cap W(\tau,U)$ may therefore be applied to this term instead.
\end{remark}
\begin{remark}\label{specialcasecalc}
    While calculating $D_S \bigcap W(\tau, U)$, if the maximal spot is $n = \sigma_t$, $\tau_{t-1} = 1$ and $t - 1 \in U$, then 
\begin{align*}
    \#\bigl(Z_m(S, \tau, U)\bigr) = \#\bigl(Z_m(S, \tau', U)\bigr) + \#\bigl(Z_m(S, \tau, U')\bigr)
\end{align*}
where $\tau'_i = \tau_i$ for $i \ne t - 1$ and $\tau'_{t-1} = 2$, and $U' = U - \{t-1\}$. This adjustment needs to be made because $\phi$ and $\psi$ are defined differently for $\tau_{t-1} = 1$ and $\tau_{t-1} \ge 2.$
\end{remark}
\begin{remark}
    While calculating $D_S \bigcap W(\tau, U)$ for $U = \emptyset$ and $\tau \in \mathbb{N}^n$, it might be that the only maximal spot at the first generation is $n$, and similarly for all the later generation until $n = s_d+\tau_{s_d}$. Therefore, $|D_S \bigcap W\big((\tau_1,\dots,\tau_{s_d})\sqcup(n-s_d)!, \phi\big)| = |D_S \bigcap W\big((\tau_1,\dots, \tau_{s_d} \sqcup (\tau_{s_d})!, \phi\big)|$. For example for $S = \{2\}, \tau_1 = 1, \tau_2=2$, we have $\bigl(1,2|(n-2)!\bigr) = (1,2|2,1)$, and the only satisfying permutation is $(1,4,2,3)$. For any $n \ge 5,$ the only satisfying permutation is $(1,4,2,3,5,\dots,n-1,n)$. In general, we can also count $|D_S \bigcap W\big((\tau_1,\dots, \tau_{s_d} \sqcup (\tau_d)!, \phi\big)|$ using our recursion formula thus get $|D_S \bigcap W\big((\tau_1,\dots,\tau_{s_d})\sqcup(n-s_d)!, \phi\big)|$ in the end.
\end{remark}
\section{Stability Range and Polynomial Degree} \label{sec:sharp}
We have established that
\begin{small}
                        \begin{align} \label{formula2a}
                    \sum\limits_{\substack{S \subseteq \{1,2,...,b\} \\ \text{ }s_1+s_2+...+s_d = b}}
   \sum_{U \subseteq \{1,2,...,s_d\}}
   \sum_{k=0}^a
   \sum\limits_{\substack{
    \tau_l = k+1 \text{ for } l \in U \\
    \tau_l \le k \text{ for } l \notin U \\
    \tau \text{ permissible}
}}
   \Bigg(\bigg|D_S \bigcap W(\tau;U)\bigg|\Bigg) \cdot 
   \Bigg([q^k] \bigg(\prod_{i=1}^{s_d} [\tau_i]_q\bigg)\Bigg)  \cdot 
   \Bigg([q^{a-k}]\bigg([n-s_d]_q!\bigg)\Bigg)
                \end{align}
            \end{small} 
            is a polynomial and the dimension of the subspaces $DR_n^{a,b}$ for big enough $n$, and in this section we will find out the earliest $n$ that the polynomial stabilizes as the dimension of the subspaces $DR_n^{a,b}$. First we prove that for $n \ge a+b$, the formula holds.
\begin{theorem}
    The polynomial starts to stabilize from $a+b$.
\end{theorem}
\begin{proof} 
    We first investigate the case when $b = 0$. Indeed, the Knuth Formula and thus Theorem~\ref{thm:knuth} is only valid for $n \ge a$, thus proving the case for $b = 0$.\\
    Now assume $b > 0$. For the first factor $|D_S\bigcap W(\tau,U)|$ in \eqref{formula2a}, the polynomial obtained through Lemma \ref{lem:telescope}, \ref{lem:powersum} and Corollary \ref{cor:polysum} is true whenever we have $n$ being a value where the recurrence has started to take place, as mentioned in the proof of Theorem \ref{maincount}. Let $s_d$ denote the biggest descent, then 
    \begin{align} \label{someeq}
        |D_S \cap W(\tau,U)| = \sum_{l = s_d + \tau_{s_d}}^n Q(l)
    \end{align}
    where $Q(n)$ is the polynomial for the sum of first-generation non-terminal nodes of the second kind. For any $s_d < b$ or $\tau_d < a+1$, equation~\eqref{someeq} yields a polynomial in $n = a+b$, since the recurrence is already valid. The only exception is the boundary case $s_d = b$ and $\tau_{s_d} = a+1$, where the recurrence has not yet become valid. This occurs when $S = \{b\}$, $b \in U$, and $k = a$, since $b \in U$ implies that $\tau_b = k+1 = a+1$ in our formula. However, in this case for any $a \ge 0, b > 0$ we have
\[
\sum_{l=a+b+1}^{a+b} Q(l) = 0
\]
by Lemma~\ref{lem:telescope}. Indeed, in this case $|D_S \cap W(\tau,U)| = 0$, so our polynomial performs the correct calculation.\\
The second factor $\Big([q^k] \big(\prod_{i=1}^{s_d} [\tau_i]_q\big)\Big)$ is a direct coefficient extraction, so it always give the correct value regardless of value of $n$. For the third parenthesis we use Knuth's formula, which gives the correct coefficient of $[q^{a-k}]$ for all $k = 0, 1, \dots, a$ values when $n \ge a+b$. Therefore our formula holds true from $n = a+b$, and stabilization starts from $a+b$ the latest.    
\end{proof}
We now know that starting from \( n = a + b \), the dimension stabilizes to our polynomial formula, but it is not yet established whether this is the earliest point at which stabilization occurs—that is, whether \( a + b \) is a sharp bound. Let $\mathscr{P}(n)$ denote the polynomial calculated using Theorem~\ref{thm:main_theorem}. To demonstrate that \( a + b \) is indeed sharp, we want to show that the polynomial $\mathscr{P}(n)$ consistently gives an incorrect value at \( n = a + b - 1 \). Specifically, we have  
\[
\mathscr{P}(a + b - 1) < \dim\left(DR_{a+b-1}^{a,b}\right).
\]

\begin{theorem}\label{sharps}
    The dimension of \( DR_n^{a,b} \) stabilizes exactly at \( a + b \); that is, \( a + b \) is a sharp bound for stabilization.
\end{theorem}

The proof of Theorem~\ref{sharps} relies on the following two lemmas.
\begin{lemma}\label{lem:knuth-conj}
    Write
    \begin{equation}\label{eq:F}
    \begin{split}
        F(n,k) = \binom{n+k-1}{k}
        &+ \sum_{j=1}^{\infty} (-1)^j \binom{n+k-u_j-1}{k-u_j} 
        + \sum_{j=1}^{\infty} (-1)^j \binom{n+k-u_j-j-1}{k-u_j-j},
    \end{split}
    \end{equation}
    where $u_j = \frac{j(3j-1)}{2}$ are the pentagonal numbers. Then
    $[q^k]\bigl([n]_q!\bigr) = F(n,k)$ for $k \le n$, while
    $F(n,n+1) = [q^{n+1}]\bigl([n]_q!\bigr) - 1$.
\end{lemma}
\begin{proof}
    We first investigate the Knuth's formula. By Euler's Pentagonal Number Theorem,
    \begin{align}
        \prod_{i=1}^\infty (1-q^k) = 1 + \sum_{k=1}^\infty (-1)^k \big(q^{\frac{3k^2-k}{2}} + q^{\frac{3k^2+k}{2}} \big)
    \end{align}
    By binomial theorem, we have that 
    \begin{equation}
        \frac{1}{(1-q)^n} = \sum_{j=0}^\infty q^j \binom{n+j-1}{j}
    \end{equation}
    By collecting the terms that have power $q^j$ we obtain Knuth's formula \cite{knuth1997art}.\\
    Notice that $[n]_q! = \frac{(1-q) \cdot (1-q^2) \dots (1-q^n)}{(1-q)^n}$, so for $m \le n$, we have that 
    \begin{align}
        [q^m] \bigg([n]_q!\bigg) &= [q^m] \Bigg( \bigg((1-q) \cdot (1-q^2) \cdot \dots \cdot (1-q^n) \bigg) \cdot \sum_{j=0}^\infty q^j \binom{n+j-1}{j}\Bigg)\\
        &= [q^m] \bigg(\prod_{i=1}^\infty (1-q^k)  \cdot \sum_{j=0}^\infty q^j \binom{n+j-1}{j} \bigg)
    \end{align}
    since any term after $(1-q^n)$ does not contribute to the coefficient of $q^m$. However for if $m = n + 1$, then Knuth's formula would give us 
    \begin{align}
         [q^m] \Bigg(\big(\prod_{i=1}^\infty (1-q^k) \big) \cdot \sum_{j=0}^\infty q^j \binom{n+j-1}{j} \Bigg) = 
        [q^m] \Bigg((1-q)  \dots (1-q^n) \cdot (1-q^{n+1}) \cdot \sum_{j=0}^\infty q^j \binom{n+j-1}{j} \Bigg) 
    \end{align}
    when we actually want 
    \begin{align}
        [q^m]\bigg([n]_q!\bigg) = [q^m] \Bigg((1-q)  \dots (1-q^n)  \cdot \sum_{j=0}^\infty q^j \binom{n+j-1}{j} \Bigg) 
    \end{align}
    Notice that the extra $(1-q^{n+1})$ contribute exactly $-1$ to the final result, so $F(n, n+1) = [q^{n+1}] \big([n]_q!\big) - 1$.
\end{proof}
\begin{lemma} \label{firstparenth}
    For the special case of $S = \{b\}, k = a,$ and $b\in U$, if we denote the polynomial we obtained for calculating $|D_S\bigcap W(\tau,U)|$ as $P(n)$, then $\big|D_S \cap W(\tau,U)\big| = 0$, but $P(a+b-1) < 0$.
   
\end{lemma}
\begin{proof}
First, we investigate the cases when $a > 0$. 
There are two possible cases for $U$ to consider: $U=\{b\}$ and
$U=\{i,i+1,\ldots,b-1,b\}$ for some $i\in\{1,\ldots,b\}$. Indeed, we do
not need to consider the case where $i\in U$ for some $i<b$ but
$U \subset \{i, i+1, \dots, b-1,b\}$. Since all positions $1,\ldots,b$ lie in the first run, we
have $\sigma_i<\sigma_{i+1}$ whenever $i<b$. Thus, if $i\in U$, then
$i+1\in U$ . Repeatedly applying this implication shows that if $i \in U$, then 
$U=\{i,i+1,\ldots,b\}$. Hence, given
$S=\{b\}$ and $b\in U$, the only possibilities are $U=\{b\}$ or
$U=\{i,i+1,\ldots,b\}$ for some $i\in\{1,\ldots,b-1\}$. We now construct
the corresponding trees in the two cases using Theorem~\ref{phimap} and \ref{psimap}.
Below is the tree for the first case:
        \begin{center}
\begin{forest}
  for tree={
    grow=south,
    s sep=50pt,
    l sep=30pt,
  }
  [{$\#\bigl(\tau_1 , \tau_2, \ldots, \ge \tau_{b-1}, \ge \tau_b \mid (n-b)!\bigr)$}
    [{$p(n) = \#\bigl(\tau_1 - 1, \tau_2 -1, \ldots, \ge \tau_{b-1} - 1 \mid (n-b)!\bigr)$}]
    [{recurrence step}]
  ]
\end{forest}
\end{center}

Below is the tree for the second case when $\tau_{b-1} > 1$:

\begin{center}
\begin{forest}
  for tree={
    grow=south,
    s sep=50pt,
    l sep=30pt,
  }
  [{$\#\bigl(\tau_1, \tau_2, \ldots, \tau_{b-1}, \ge \tau_b \mid (n-b)!\bigr)$}
    [{$p(n) = \#\bigl(\tau_1-1, \tau_2-1, \ldots, \tau_{b-1} -1\mid (n-b)!\bigr)$}]
    [{recurrence step}]
  ]
\end{forest}
\end{center}
Below is the tree for the second case when $\tau_{b-1} = 1$. Note permissibility would force $\tau_{b-2} = 2, \dots, \tau_1 = b-1$.
\begin{center}
\begin{forest}
  for tree={
    grow=south,
    s sep=50pt,
    l sep=30pt,
  }
  [{$\#\bigl(b-1, b-2, \ldots,2, 1, \ge \tau_b \mid (n-b)!\bigr)$}
    [{$p(n) = \#\bigl(n-1, n-2, \dots, 2,1\bigr)$}]
    [{recurrence step}]
  ]
\end{forest}
\end{center}
    For $D_S \cap W(\tau,U)$ to be non-empty, we need the length of the last run $n-b$ to be at least $\tau_b = a+1$ which is equivalent to $n \ge a+b+1$, and we know that for $n = a+b$, we have $\big|D_S \cap W(\tau,U)\big| = 0$. Now we prove that our polynomial however, gives a negative number for $n = a+b-1$, i.e. $P(a+b-1) < 0$.\\
    Before we investigate case by case, notice that we only have one maximal spot $b$ thus one non-recurrence child being $Z_b(\{b\}, \tau, U)$, so if $p(n)$ counts the non-terminal node of the second kind child, we have that 
    \begin{align}\label{nab}
        P(n) =  \sum_{i= a+b+1}^n p(i) &= \sum_{i=1}^{n} p(i) - \sum_{i=1}^{a+b} p(i)
    \end{align}
    so when $n = a+b-1$, we have that $P(n) = -p(a+b)$. Now we investigate the two cases mentioned above.
    \begin{enumerate}
        \item For this case, we have $D_{\{b\}} \cap W(\tau, U)$, and $U \supset \{b\}$, so $\tau_{b-1} = a+1$. We know that the non-terminal node of the second kind child at the first generation is
        \begin{align}
            Z_b\big(\{b\}, \tau, U) = \bigl(\tau_1, \tau_2, \ldots, \ge \tau_{b-1}-1 \mid (n-b)!\bigr)
        \end{align}
        Let $p(n) = \#\bigl(\tau_1, \tau_2, \ldots, \ge \tau_{b-1}-1 \mid (n-b)!\bigr)$, then the smallest length of the last run $(n-b)!$ for $p(n) > 0$ will be the $\tau_{b-1}-1$. In our case $p(n) > 0$ when $n - b \ge \tau_{b-1}-1$ which is equivalent to $n \ge a + b$. Therefore $p(a+b) > 0,$ so $P(a+b-1) < 0$.
        \item For this case, we have $D_{\{b\}} \bigcap W(\tau, \{b\})$, and $\tau_{b-1} \le a$. We first investigate when $\tau_{b-1} > 1$. We know that the non-terminal node of the second kind child at the first generation is
        \begin{align}
            Z_b\big(\{b\}, \tau, \{b\}) = \bigl(\tau_1, \tau_2, \ldots, \tau_{b-1} -1\mid (n-b)!\bigr)
        \end{align}
        Let $p(n) = \#\bigl(\tau_1, \tau_2, \ldots, \tau_{b-1} -1 \mid (n-b)!\bigr)$. \\
        If $\tau_{b-1} > 1,$ then with same reasoning above we know $p(n) > 0 $ when the length of the last run $n-b \ge \tau_{b-1} -1$, which is equivalent to $n \ge \tau_{b-1}+b-1$. Since $\tau_{b-1} \le a,$ in particular we have $p(a+b) > 0,$ so $P(a+b-1) <0.$ \\
        Now assume $\tau_{b-1}=1$. The non-terminal node of the second kind child at the first generation consists only of the identity permutation, obtained by merging the first and last runs after removing $n$ at the end of the first run. Thus,
  $p(n) = 1 > 0$ for all $n\ge2$. Since $b-1\notin U$ and $\tau_{b-1}\le k\le a$, we have $a\ge1$. Together with our assumption that $b\ge1$, this gives $a+b\ge2$, and hence $p(a+b)>0$. Therefore,
 $P(a+b-1) < 0$ in this case too.
    \end{enumerate}
    Now we assume $a = 0.$ Since $\tau_i = 0$ if $i \not\in U$, we know that the only nonzero $\bigl|D_S \cap W(\tau,U)\bigr|$ is when $U = \{1, 2,\dots, b\}$, so we are calculating $\#\bigl(\ge1, \ge1,\dots, \ge1|(n-b)!\bigr)$. By Remark~\ref{permissvsU}, this is equivalent to counting
\[
P(n)=\#\bigl(\ge b-1,\ge b-2,\dots,\ge2,\ge1,\ge1|(n-b)!\bigr),
\]
which is precisely the number of permutations whose descent set is $\{b\}$.
We can count these permutations directly. Such a permutation $\sigma$ is determined by choosing the set $\{\sigma_1,\dots,\sigma_b\}$ and arranging its elements in increasing order in positions $1,2,\dots, b$, while arranging the remaining elements in increasing order in positions $b+1,\dots,n$. Hence, there are $\binom{n}{b}$ possible permutations. This construction guarantees that the descent set is contained in $\{b\}$, but it does not guarantee that $b$ is indeed a descent; that is, we still need to impose $\sigma_b>\sigma_{b+1}$. If instead $\sigma_b<\sigma_{b+1}$, then
\(
\sigma_1<\cdots<\sigma_b<\sigma_{b+1}<\cdots<\sigma_n,
\)
so $\sigma$ must be the identity permutation. Thus, exactly one of the $\binom{n}{b}$ permutations must be excluded, and consequently
\(
P(n)=\binom{n}{b}-1.
\)
In particular, when $a=0$ and $n=a+b-1=b-1$, we obtain
\(
P(a+b-1)=\binom{b-1}{b}-1=-1<0.
\)
\end{proof}
Now we are ready to prove Theorem~\ref{sharps}. 
\begin{proof}
We first investigate the cases when $b = 0$.
For $b = 0,$ we use \eqref{b0formula}, and denote the formula as $\mathscr Q(n).$ By Lemma~\ref{lem:knuth-conj}, we know that 
\[ \mathscr{Q}(a-1)  = F(a-1,a) < [q^{a}]\bigl([a-1]_q!\bigr) = dim(DR_{a-1}^{a,b}).\]
Since $\mathscr{Q}(n)=\dim DR_n^{a,0}$ for all $n\ge a$, any polynomial that agrees with $\dim DR_n^{a,0}$ for all $n\ge a-1$ must agree with $\mathscr{Q}(n)$ for infinitely many values of $n$, and hence must equal $\mathscr{Q}(n)$. Since
\(
\mathscr{Q}(a-1)<\dim DR_{a-1}^{a,0},
\)
this is impossible. Thus $a$ is the sharp stabilization bound.
Now we turn to the case $b>0$. We first consider the case $a>0$.\\
For each choice of
$(S,U,\tau,k)$, denote the corresponding summand in \eqref{formula2a} by
\[
C_{S,U,\tau,k}(n)
=
\big|D_S\cap W(\tau,U)\big|
\cdot
[q^k]\!\left(\prod_{i=1}^{s_d}[\tau_i]_q\right)
\cdot
[q^{a-k}]\bigl([n-s_d]_q!\bigr).
\]
The middle factor is always computed exactly by coefficient extraction, so
possible discrepancies at $n=a+b-1$ can only arise from the first factor
$\big|D_S\cap W(\tau,U)\big|$ or the third factor
$[q^{a-k}]([n-s_d]_q!)$ in the formula. We analyze these two possibilities separately.\\
We first consider the boundary cases of the first type in which the calculation of the first factor fails to
be correct. Notice that, for $n<a+b$, the algorithm of
Theorem~\ref{maincount} gives the correct polynomial for
$\big|D_S\cap W(\tau,U)\big|$ for every choice of
$(S,U,\tau,k)$ except when
\[
S=\{b\},\qquad k=a,\qquad b\in U.
\]
This is because for $n \ge \tau_{s_d} + s_d = a+b+1$, the recurrence has set in so the polynomial is positive and valid; for $n = \tau_{s_d} + s_d -1 = a+b,$ we have set it to be 0, which is the correct value for counting $|D_S\cap W(\tau,U)|.$ Therefore as long as $n \ge \tau_{s_d} + s_d - 1$, we have $P(n) = |D_S\cap W(\tau,U)|$. For example, if $b=3$, $S=\{3\}$, and $k=a=2$, then the polynomial
obtained for $11\ge 3|(n-3)!$ is valid for $n\ge6$, and remains valid for
$n=5$ by Lemma~\ref{lem:telescope}. However, by
Lemma~\ref{firstparenth}, it becomes negative for $n=4$, whereas the true
value of $\big|D_S\cap W(\tau,U)\big|$ is $0$.\\
Moreover, this boundary is sharp. If we instead take $k=1$, then the
algorithm still produces the correct polynomial at $n=4=a+b-1$. Likewise,
if we replace $S=\{3\}$ by $S=\{1,2\}$, then $s_d=2$, so
$\tau_{s_d}+s_d=4$, and Lemma~\ref{lem:telescope} again guarantees that the
polynomial remains valid at $n=4$. Therefore, the only boundary cases of
this type are those considered in Lemma~\ref{firstparenth}.\\
For each such boundary case, we know that the true value is
\(
\big|D_S\cap W(\tau,U)\big|=0,
\)
while the polynomial $P(n)$ produced by Theorem~\ref{maincount} is negative at
$n=a+b-1$. The remaining two factors of
$C_{S,U,\tau,k}(a+b-1)$ are nevertheless computed correctly and they are nonzero: the second
factor is nonzero since $\tau_b=a+1$, and the third factor is the
constant term of $[n-b]_q!$, which is $1$. Hence every such summand
$C_{S,U,\tau,k}(a+b-1)$ is strictly smaller than its true contribution.\\
Next we consider the boundary cases of the second type, arising from the third factor. When
$n=a+b-1$, this occurs only when
\[
S=\{b\},\qquad k=0, \qquad U \subseteq \{1,2,\dots,b-1,b\}
\]
so that the third factor is
\(
[q^a]([a-1]_q!).
\)
By Lemma~\ref{lem:knuth-conj}, Knuth's formula underestimates this coefficient.
Since these are not boundary cases of the first type for $a > 0$, the first factor is
computed correctly by Theorem~\ref{maincount}, while the second factor is
again computed exactly. Therefore for this type of boundary cases, every corresponding summand
$C_{S,U,\tau,k}(a+b-1)$ is also strictly smaller than its true value.\\
Finally, every non-boundary summand of this type agrees with its true value at
$n=a+b-1$, whereas every boundary summand of this type contributes strictly less than its
true value. Consequently, for $a>0,$ we have
\(
\mathscr P(a+b-1)
<
\dim DR_{a+b-1}^{a,b}.
\)\\
Now suppose that $a=0$. The coefficient of $q^0t^b$ in the Schedules formula \eqref{eq:schedules} is the number of permutations $\sigma\in S_n$ with $\operatorname{maj}(\sigma)=b$. Therefore for our formula to be correct, the product of the second and third factors must be $1$. Since both factors are positive integers, this can occur only if each of them is equal to $1$. We consider the contributions corresponding to the possible descent sets $S$ separately.
First, suppose that $S\neq\{b\}$. In this case, the first factor in $C_{S,U,\tau,k}(b-1)$ is not a boundary case of the first type, and hence is computed correctly by our formula. Since the constant term of any product of $q$-integers is 1, the second factor $[q^0]\bigl(\prod\limits_{i=1}^{s_d} [w_i(\sigma)]_q\bigr)$ is $1$; the third factor is
\(
F(-1,0)=1.
\)
Both of the second and third factor are correct, so $C_{S,U,\tau,k}(b-1)$ gives the correct contribution for every $S\neq\{b\}$.
It remains to consider $S=\{b\}$. In this case, the first factor $\bigl(|D_S \cap W(\tau,U)\bigr|$ is a boundary case of the first type, and by Lemma~\ref{firstparenth} its polynomial expression gives
\(
P(b-1)
<
\bigl|D_S\cap W(\tau,U)\bigr|.
\)
The second factor is again $1$, and the third factor is $F(-1,0)=1$, which are both correct. Hence, the contribution
\(
C_{S,U,\tau,k}(b-1)
\)
is strictly underestimated by the polynomial formula for boundary case of this type.
There is no boundary case of the second type when $a=0$: since $k=a=0$, the third factor is 
\(
F(-1,0)=1,
\)
which is correct. Therefore, the only discrepancy occurs in the boundary case of the first type with $S=\{b\}$, and consequently
\(
\mathscr P(b-1)
<
\dim DR_{b-1}^{0,b}.
\)
Combining this with the case $a>0$, we obtain, for all $a\in\mathbb{N}\cup\{0\}, b\in \N$,
\[
\mathscr P(a+b-1)
<
\dim DR_{a+b-1}^{a,b}.
\]
Since $\mathscr{P}(n)=\dim DR_n^{a,b}$ for all $n\ge a+b$, any polynomial that agrees with $\dim DR_n^{a,b}$ for all $n\ge a+b-1$ must agree with $\mathscr{P}(n)$ for infinitely many values of $n$, and hence must equal $\mathscr{P}(n)$. Since
\(
\mathscr{P}(a+b-1)<\dim DR_{a+b-1}^{a,b},
\)
this is impossible. Thus $a+b$ is the sharp stabilization bound.
\end{proof}
\begin{remark}
    In their foundational paper, Church and Farb \cite{Church_2015} developed the theory of FI-modules and used these structures to establish results in representation stability. In particular, they introduced the notion of FI\#-modules, a more rigid structure than FI-modules. The existence of a sharp bound for the dimension polynomial demonstrates that $DR_n^{a,b}$ is not an FI\#-module—without requiring an explicit analysis of its FI-module structure.
\end{remark}
Now we investigate the degree of the polynomial once it stabilizes.
\begin{theorem}\label{thm:main}
Fix $a, b \in \mathbb{Z}_{\ge0}$. For $n \ge a+b$, $\dim(DR_n^{a,b})$is a polynomial in $n$ of degree
exactly $a + b$.
\end{theorem}

\begin{proof}
For $b = 0$, the highest degree term in \eqref{b0formula} is $\binom{n+a-1}{a}$ which is of degree $a$, so indeed the coefficient of $q^at^0$ in $Hilb(DR_n^{a,0})$is $a$. Now we focus on $b \ge 1$ when we use \eqref{formula2a}.\\
The second factor is a constant arising from coefficient extraction and
therefore does not contribute to the degree in $n$. We thus focus on the
first and third factors and determine when they simultaneously attain
their maximal degrees.\\
For the first factor $|D_S \cap W(\tau,U)|$, Theorem~\ref{maincount} shows
that at most $s_d-1$ applications of $\phi$ are required to reach the
base case $S=\{1\}$. By Corollary~\ref{cor:polysum}, if exactly
$s_d-1$ applications of $\phi$ are required, then the degree of the
resulting polynomial is
\(
(s_d-1) + \deg(\text{base case}).
\)
By Lemma~\ref{lem:base}, the base case polynomial is either degree 1 or 0. Since $\operatorname{maj}(S)=b$, we have $s_d\leq b$,
with equality precisely when $S=\{b\}$. Thus
\[
\deg\left|(D_S \cap W(\tau,U)\right|
\leq (b-1)+1=b.
\]
This maximum is attained when
\(
U=\{1,2,\ldots,b-1,b\}.
\)
Indeed, in this case, every application of $\phi$ decreases the maximum
element of the descents set by exactly one:
\(
\max \left(\phi^{i+1}(S)\right)
=
\max \left(\phi^i(S)\right)-1
\)
by Definition~\ref{phimap} and Remark~\ref{specialcasecalc}. Hence, when
$S=\{b\}$, exactly $b-1$ applications of $\phi$ are required to reach
the base case, and we have $deg\bigl(\bigl|D_S \cap W(\tau,U)\bigr|\bigr) = b$.\\
For the third factor
\(
[q^{a-k}]\left([n-s_d]_q!\right),
\)
Knuth's formula~\eqref{eq:knuth} shows that
\(
\deg\left([q^r]\left([n]_q!\right)\right)=r.
\)
Therefore,
\[
\deg\left([q^{a-k}]\left([n-s_d]_q!\right)\right)=a-k.
\]
Thus, the maximal possible degree of the third factor is $a$, achieved
when $k=0$.\\
Consequently, the maximal degree of the product of the first and third
factors is attained when
\[
S=\{b\}, \qquad
U=\{1,2,\ldots,b-1,b\}, \qquad
k=0, \qquad
\tau=
\left(
\underbrace{1,\ldots,1}_{b\text{ times}},
(n-b)!
\right).
\]
The corresponding configuration is not necessarily permissible. By the
dissection procedure described in Remark~\ref{permissvsU}, we may
replace it by the permissible configuration
\(
\tau=
\bigl(b-1,b-2,\ldots,2,1,1,(n-b)!\bigr).
\)
We have counted this in the proof of Lemma~\ref{firstparenth} to be $\binom{n}{b} - 1$, so $\bigl|D_S\cap W(\tau,U)\bigr|$ is a polynomial in $n$ of degree $b$.
When $a>0$, since $k=0$, the third factor has degree $a$. For when $a = 0,$ the third factor is $F(-1,0) = 1$ which is also of degree $a = 0$. Hence for both cases, the resulting
summand has degree $a+b$, and the dimension polynomial
$dim(DR_n^{a,b})$ has degree exactly $a+b$.
\end{proof}
\begin{remark}
    Since $DH_n \cong DR_n$ as $S_n$-modules, every result we have proved in this paper is also true for $DH_n$.
\end{remark}
We actually improved both the stability range and the polynomial degree. Bahran \cite{BAHRAN2024429} found a stability bound for \( \dim(DR_n^{a,b}) \) to be \(2(a+b)-1 \), and we found a sharp bound of \( a+b \). By using representation-theoretic tools, Church et al.\ \cite{Church_2015} found the degree of the dimension polynomial \( \dim(DR_n) \) to be at most \( a+b \), while we found it to be \textit{exactly} \( a+b \). We provide a table of polynomials calculated by using the algorithm developed and proven in this article of \( \dim(DR_n^{a,b}) = [q^a t^b] \left( \sum\limits_{\sigma \in S_n} t^{\mathrm{maj}(\sigma)} \prod\limits_{i=1}^n [w_i(\sigma)]_q \right) \), for various $a,b$ and $n \ge a+b$ (Table~1).

\begin{table}[ht]
    \centering
    \small 
    \renewcommand{\arraystretch}{1.5} 
    \begin{tabular}{|c|c|c|c|c|} 
        \hline
        \textbf{/} & \textbf{q$^0$} & \textbf{q$^1$} & \textbf{q$^2$} & \textbf{q$^3$} \\ 
        \hline
        $t^0$ & 1 & $n-1$ & $\frac{1}{2} n^2 -\frac{1}{2}n - 1$
 & $\frac{1}{6} n^3 - \frac{7}{6} n$
\\ 
        \hline
        $t^1$ & $n - 1$ & $n^2 - 2n$ & $\frac{n^3}{2} - n^2 - \frac{3n}{2} + 1$ & $\frac{1}{6}n^4 - \frac{1}{6}n^3 - \frac{5}{3}n^2 + \frac{2}{3}n + 1$ \\ 
        \hline
        $t^2$ & $\frac{n^2}{2} - \frac{n}{2} - 1$ & $\frac{n^3}{2} - n^2 - \frac{3n}{2} + 1$ & $\frac{n^4}{4} - \frac{n^3}{2} - \frac{7n^2}{4} + n + 1$ & $\frac{1}{12}n^5 - \frac{1}{12}n^4 - \frac{5}{4}n^3 + \frac{1}{12}n^2 + \frac{13}{6}n + 1$ \\ 
        \hline
        $t^3$ & $\frac{n^3}{6} - \frac{7n}{6}$ & $\frac{n^4}{6} - \frac{n^3}{6} - \frac{5n^2}{3} + \frac{2n}{3} + 1$ & $\frac{n^5}{12} - \frac{n^4}{12} - \frac{5n^3}{4} + \frac{n^2}{12} + \frac{13n}{6} + 1$ & $\frac{n^6}{36} - \frac{23n^4}{36} - \frac{n^3}{2} + \frac{19n^2}{9} + 3n - 1$ \\ 
        \hline
    \end{tabular}
    \caption{Polynomials for $dim(DR_n^{a,b})$}
    \label{tab:expressions}
\end{table}

\bibliographystyle{amsalpha}
    \bibliography{bib}

\end{document}